\pdfoutput=1
\documentclass[a4paper,12pt]{article}
\usepackage[left=2.5cm, right=2.5cm, top=3.0cm, bottom=3.0cm]{geometry}

\usepackage{lipsum}
\usepackage{amsfonts}
\usepackage{graphicx}
\usepackage{epstopdf}
\ifpdf
  \DeclareGraphicsExtensions{.eps,.pdf,.png,.jpg}
\else
  \DeclareGraphicsExtensions{.eps}
\fi

\usepackage{amsthm}
\usepackage{fancyhdr}
\usepackage{hyperref}

\usepackage{amsmath}
\usepackage{datetime}
\newdate{date}{15}{06}{2022}
\date{\showdowfalse\displaydate{date}}
\usepackage{pgfplots} 

\makeatletter
\newcommand{\Eqref}[1]{\textup{\tagform@{\ref*{#1}}}}
\makeatother

\theoremstyle{remark}

\usepackage{color}
\usepackage{xspace}
\usepackage{url}
\usepackage{graphicx}
\usepackage{lipsum} 
\usepackage{enumitem} 
\usepackage{upgreek}
\usepackage{subfigure}
\usepackage{algpseudocode}
\usepackage{longtable}
\usepackage{tabu}
\usepackage{bm}
\usepackage{adjustbox}
\usepackage{bm}
\usepackage{placeins}
\usepackage{float}
\usepackage[normalem]{ulem}
\usepackage[subtle]{savetrees}
\usepackage{booktabs}
\usepackage{pgfplots}
\usepackage{pgf-pie}

\usepackage{arydshln}

\usepackage{setspace}

\usepackage{listings}
\usepackage{textcomp}
\definecolor{lbcolor}{rgb}{0.93,0.93,0.93}
\newcommand{\PreserveBackslash}[1]{\let\temp=\\#1\let\\=\temp}
\newcolumntype{C}[1]{>{\PreserveBackslash\centering}p{#1}}

\definecolor{gnuplot@orange}{RGB}{229,158,0}
\definecolor{gnuplot@purple}{RGB}{148,0,212}
\definecolor{gnuplot@red}{RGB}{200,0,0}
\definecolor{gnuplot@lightblue}{RGB}{87,181,232}
\definecolor{gnuplot@green}{RGB}{0,158,65}
\definecolor{gnuplot@darkblue}{RGB}{0,115,179}
\definecolor{gnuplot@yellow}{RGB}{240,227,66}
\pgfplotscreateplotcyclelist{colorGPL}{%
gnuplot@darkblue,every mark/.append style={fill=gnuplot@darkblue!80!black},mark=o\\%
gnuplot@red,every mark/.append style={fill=gnuplot@red!50!black},mark=square\\%
gnuplot@green,every mark/.append style={fill=gnuplot@green!50!black},mark=triangle*\\%
gnuplot@orange,every mark/.append style={fill=gnuplot@orange!80!black,mark size=2.5pt},mark=x\\%
gnuplot@purple,mark=diamond\\%
black,densely dashed,every mark/.append style={solid,fill=gnuplot@darkblue!80!black},mark=square*\\%
gnuplot@lightblue,densely dashed,every mark/.append style={solid,fill=gnuplot@lightblue!30!black},mark=otimes*\\%
red!80!white,densely dashed,every mark/.append style={solid,fill=red!40!white},mark=oplus*\\%
}

\pgfplotsset{compat=1.9}

\newcommand{\myred}{black}

\definecolor{brightmaroon}{rgb}{0.76, 0.13, 0.28}
\definecolor{dukeblue}{rgb}{0.0, 0.0, 0.61}

\graphicspath{{svg_with_tex/}}

\title{Stage-parallel fully implicit Runge-Kutta implementations with optimal multilevel preconditioners at the scaling limit\thanks{Dedicated
to the memory of Owe Axelsson.}}

\author{Peter Munch\thanks{Corresponding author. Helmholtz-Zentrum Hereon, University of Augsburg (peter.muench@uni-a.de).}
\and Ivo Dravins\thanks{Uppsala University (ivo.dravins/maya.neytcheva@it.uu.se).}
\and Martin Kronbichler\thanks{University of Augsburg, Uppsala University (martin.kronbichler@uni-a.de).}
\and Maya Neytcheva\footnotemark[3]
}
  
\usepackage{amsopn}

\ifpdf
\hypersetup{
  pdftitle={Stage-parallel fully implicit Runge-Kutta implementations with optimal multilevel preconditioners at the scaling limit},
  pdfauthor={P. Munch, I. Dravins, M. Kronbichler, and M. Neytcheva}
}
\fi

\begin{document}

\maketitle

\begin{abstract}
We present an implementation of a fully stage-parallel preconditioner
for Radau IIA type {\color{\myred}fully implicit Runge--Kutta methods, which approximates
the inverse of
$A_Q$ from the Butcher tableau by the lower triangular matrix resulting from an LU decomposition and diagonalizes the system with as many blocks as stages. For the transformed system, we employ} a block preconditioner where
each block is {distributed and} solved by a subgroup of processes in parallel. For combination
of partial results, we either use a communication pattern resembling Cannon's algorithm or shared memory.
A {performance model and a} large set of performance studies (including strong scaling runs with up to
150k processes on 3k compute nodes) conducted for a time-dependent heat problem,
using matrix-free finite element methods, indicate that the
stage-parallel implementation can reach higher throughputs when the block solvers operate at lower parallel
efficiencies, which occurs near the scaling limit. Achievable speedup increases linearly
with number of stages {and are bounded by the number of stages}. {\color{\myred}Furthermore, we  show that
the presented stage-parallel concepts are also applicable to the case  that
 $A_Q$ is directly diagonalized, which requires complex arithmetic or 
{the solution of two-by-two blocks}
and sequentializes parts of the algorithm.} {Alternatively to
distributing stages and assigning them to distinct processes, we discuss
the possibility of batching operations from different stages together.}
\end{abstract}



\section{Introduction}

Runge--Kutta methods are widely used time-integration
schemes to solve ordinary differential equations (ODE) of the form:
\begin{align*}
    \frac{dy}{dt} = f(t,y).
\end{align*}
We consider a partial differential equation, rewritten as a system of ODEs after using finite element methods (FEM) to discretize in space.
The basic algorithm is to obtain the solution at the next time step
via linear combination of $Q$ intermediate stage solutions:
\begin{align*}
\textbf{y}_{n+1} = \textbf{y}_n + \tau \sum_{1 \le i \le Q} b_i \textbf{k}_i
\quad
\text{with}
\quad
\textbf{k}_i = f( t_n + c_i \tau , \textbf{y}_n + \tau \sum_{1 \le j \le Q} a_{ij} \textbf{k}_j ),
\end{align*}
where $t_n$ is the time at time step $n$ and $\tau$ the current time-step size.
The Butcher tableau is a compact notation for these methods in terms of a matrix $A_Q$ as well as two vectors $\bm{b}_Q$ and $\bm{c}_Q$.

There is a vast literature on optimizing Runge--Kutta methods. The investigations include improving
the accuracy and the stability region as well as performance optimization. Low-storage Runge--Kutta methods~\cite{Kennedy00}, for instance,
update the solution step by step so that the intermediate results for all stages need not be stored simultaneously, which
might be an advantage for memory-intensive applications, e.g., computational plasma physics \cite{munch2021hyper}.
{\color{\myred} Implicit Runge--Kutta methods
solve systems of increased size so that the development of {efficient solvers} is crucial~\cite{abu2022monolithic, Axelsson1, AxelDravNeyt, bickart1977efficient, burrage1999parallel, butcher1976implementation, de1998specification, jay1999parallelizable, pazner2017stage, southworth2022fast2, southworth2022fast}.}

\begin{figure}[!t]

\centering

 \begin{tikzpicture}[thick,scale=0.9, every node/.style={scale=0.75}]
    \begin{loglogaxis}[
      width=0.48\textwidth,
      height=0.32\textwidth,
      title style={font=\tiny},every axis title/.style={above left,at={(1,1)},draw=black,fill=white},
      xlabel={Nodes ($\times$ 48 CPU cores)},
      ylabel={Time per time step [s]},
      legend pos={south west},
      legend cell align={left},
      cycle list name=colorGPL,
      grid,
      semithick,
      ymin=0.0008,ymax=10,
      xmin=0.8,xmax=4096,
      xtick={1, 2,4, 8, 16,32,64,128,256,512, 1024, 2048, 3072},
      xticklabels={1, ,4, , 16,,64,,256,, 1k, ,3k},
	mark options={solid}
      ]


\addplot+[solid, gnuplot@red,every mark/.append style={fill=gnuplot@red!80!black},mark=square] table [x index=0,y index=1] {data/gmg/data_0007_0.tex};

\addplot+[solid, gnuplot@green,every mark/.append style={fill=gnuplot@green!80!black},mark=square] table [x index=0,y index=1] {data/gmg/data_0008_0.tex};

\addplot+[solid, magenta,every mark/.append style={fill=magenta!80!black},mark=square] table [x index=0,y index=1] {data/gmg/data_0009_0.tex};

\addplot+[solid, brown,every mark/.append style={fill=magenta!80!black},mark=square] table [x index=0,y index=1] {data/gmg/data_0010_0.tex};

\addplot+[solid, blue,every mark/.append style={fill=magenta!80!black},mark=square] table [x index=0,y index=1] {data/gmg/data_0011_0.tex};

\addplot+[solid, orange,every mark/.append style={fill=magenta!80!black},mark=square] table [x index=0,y index=1] {data/gmg/data_0012_0.tex};

    \end{loglogaxis}
\end{tikzpicture}\qquad
\begin{tikzpicture}[thick,scale=0.9, every node/.style={scale=0.75}]
    \begin{loglogaxis}[
      width=0.48\textwidth,
      height=0.32\textwidth,
      title style={font=\tiny},every axis title/.style={above left,at={(1,1)},draw=black,fill=white},
      xlabel={Nodes ($\times$ 48 CPU cores)},
      ylabel={Time per time step [s]},
      legend pos={south west},
      legend cell align={left},
      cycle list name=colorGPL,
      grid,
      semithick,
      ymin=0.0008,ymax=10,
      xmin=0.8,xmax=4096,
      xtick={1, 2,4, 8, 16,32,64,128,256,512, 1024, 2048, 3072},
      xticklabels={1, ,4, , 16,,64,,256,, 1k, ,3k},
	mark options={solid}, legend pos=north east
      ]

%
%

\addplot+[solid, magenta,every mark/.append style={fill=magenta!80!black},mark=square] table [x index=0,y index=1] {data/gmg/data_0009_0.tex};
\addplot+[dashed, magenta,every mark/.append style={fill=magenta!80!black},mark=x] table [x index=0,y expr= \thisrow{time}*8] {data/gmg/data_0009_0.tex};

%
%


\draw[draw=black,thick] (axis cs:64, 7.208600e-03) -- node[above,rotate=-0,inner sep=2pt,outer sep=0.5pt]{}(axis cs:2048, 7.208600e-03);
\draw[draw=black,thick] (axis cs:512, 0.0287072) -- node[above,rotate=-0,inner sep=2pt,outer sep=0.5pt]{}(axis cs:2048, 0.0287072);
\draw[draw=black,thick,<->] (axis cs:1024, 7.208600e-03) -- node[above,rotate=-0,inner sep=2pt,outer sep=0.5pt]{}(axis cs:1024, 0.0287072);
\node (a) at (axis cs:1800,0.015) {$\times 3.98$};

\draw[draw=black,thick] (axis cs:8, 4.050200e-02) -- node[above,rotate=-0,inner sep=2pt,outer sep=0.5pt]{}(axis cs:160, 4.050200e-02);
\draw[draw=black,thick] (axis cs:64, 0.0576688) -- node[above,rotate=-0,inner sep=2pt,outer sep=0.5pt]{}(axis cs:160, 0.0576688);
\draw[draw=black,thick,<->] (axis cs:128, 4.050200e-02) -- node[above,rotate=-0,inner sep=2pt,outer sep=0.5pt]{}(axis cs:128, 0.0576688);
\node (a) at (axis cs:270,0.048) {$\times 1.42$};

\draw[draw=black,thick] (axis cs:1, 3.050300e-01) -- node[above,rotate=-0,inner sep=2pt,outer sep=0.5pt]{}(axis cs:16, 3.050300e-01);
\draw[draw=black,thick] (axis cs:8, 0.324016) -- node[above,rotate=-0,inner sep=2pt,outer sep=0.5pt]{}(axis cs:16, 0.324016);
\node (a) at (axis cs:32,0.32) {$\times 1.06$};

\legend{1 iteration, 8 iterations}

    \end{loglogaxis}
\end{tikzpicture}

\caption{a) Time of a single conjugate gradient iteration preconditioned by GMG for different numbers of refinements $6 \le L \le 11$ and $k=1$. b) Example visualizing the benefit of stage parallelism for $Q=8$, by comparing the
time of a single iteration with the one of eight iterations and providing idealized speedups
in the case that $Q$ iterations are solved in parallel by $Q$ subgroups.}\label{fig:intro:gmg}

\end{figure}
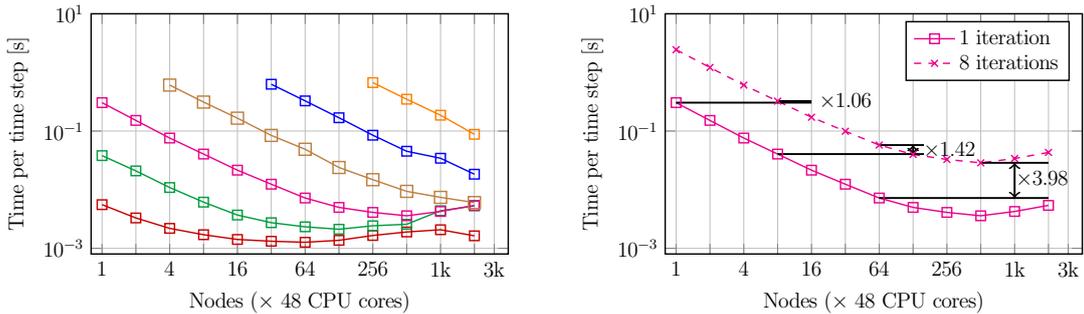

In this work, {\color{\myred}we investigate {solvers} for fully implicit Runge--Kutta methods, focusing on the parallelization of the solution process
of each stage \cite{JacksonNorsett}}. Stage-parallel approaches and, similarly, parallel-in-time approaches~\cite{Bolten, Lions} (out of the scope of this work) are motivated by
the scaling limit of algorithms, e.g., of iterative solvers, on distributed systems when parallelism is only exploited in the spatial domain. In the latter case, algorithms can not
be run faster than a certain threshold even if more hardware resources are added.
Figure~\ref{fig:intro:gmg}a) shows, as an example, the times of one conjugate gradient iteration
preconditioned by a single V-cycle of geometric multigrid (GMG)
we will use in the experimental section
when run across a large
range of number of processes.
It is clearly visible that the times flatten out when the work per process
becomes too little. This phenomenon is also
known as ``scaling limit''. In the context of GMG, the minimum time of one iteration is approximately proportional
to the number of levels. {Since the total solution time of IRK is approximately
the solution time accumulated over all stages
run sequentially, the behavior of the overall solution time is similar under the
assumption that each block corresponding to a stage can be solved by GMG.}

As an alternative to the sequential solution process of each stage, one could 
assume that it is possible to solve the stages by $Q$ process groups with GMG independently.
In consequence, while the time for solving for one stage might take longer, the stages
can be solved in parallel. A speedup can be expected if the
solution time of one stage does not increase significantly ($\ll Q$). 
Figure~\ref{fig:intro:gmg}b) presents, as an example, the
solution times for a single GMG iteration and for eight GMG iterations (assuming $Q=8$).
In addition, it shows ideal speedups we can expect by solving one
GMG iteration in parallel by
each subgroup with one eighth of the processes. The
speedup increases with increasing number of processes. Away from the scaling
limit, the ideal speedup is only minor (a few percent), which will be dominated by organizational overheads in practice. At the scaling limit, the ideal speedup
approaches $Q$.
In summary, {in the case that the IRK algorithm could be reformulated
such that stages can be solved independently,} this additional level of parallelism might allow to
increase the granularity of the subproblem to solve and better
utilize the capacity of the given hardware resources, such as to reach lower times to solution. {In this publication, we show that such a reformulation is possible and, for a simple
benchmark, we can achieve significant speedup at the scaling limit this way.}

{\color{\myred}Stage-parallel implementations of implicit Runge--Kutta methods have been
investigated but  {\color{\myred}rarely} implemented in the literature. One example is the work by Pazner and Persson~\cite{pazner2017stage}, who
considered a stage-parallel IRK preconditioner using a block-Jacobi solver
across the processes around a local ILU: by solving the stages by subgroups of processes,
the size of each block to be solved was increased and the efficiency of ILU was
improved in terms of number of iterations so that speedups were reached.
However, a critical discussion on the benefits of stage-parallel IRK regarding performance and
on the challenges regarding efficient implementation
in the context of more sophisticated global preconditioners and optimal preconditioners
for the blocks, such as multigrid,
are lacking in the literature.\footnote{We define ``optimal'' as solver  with high node-level performance whose number of iterations
is independent of the number of DoFs and of the number of processes.}}
In order to address this issue, {we consider a direct
factorization of the linear system arising from the implicit Runge--Kutta method and have extended a novel preconditioner for IRK},
which has been introduced in \cite{axelsson2020numerical}.
Our results are based on benchmark programs
leveraging the infrastructure of the open-source FEM library \texttt{deal.II}~\cite{dealII93, dealii2019design}
and are available on GitHub at
\url{https://github.com/peterrum/dealii-spirk}.

{As a critical remark, we note that running IRK and, generally, time-stepping schemes
at the scaling limit means that the user has enough hardware resources and computational
budget. Most scientists, however, are not in such a ``luxurious'' position. Hence, using
stage-parallel Runge-Kutta methods and potentially parallel-in-time algorithms may not benefit them or might even be
disadvantageous if applied far from the scaling limit due to some---unavoidable---organizational overhead.}

The remainder of this contribution is organized as follows. Section~\ref{sec:methods} provides a short
description of the linear systems arising from the implicit Runge--Kutta method and {of stage-parallel solution procedures, followed by a discussion of their building blocks}. Section~\ref{sec:implementation} presents
implementation details of the building blocks, and Section~\ref{sec:modeling} discusses relevant performance models.
{Sections~\ref{sec:performance}--\ref{sec:complex} demonstrate performance results for the solution of a time-dependent
heat equation discretized by low- and high-order FEM with different stage-parallel
solution procedures and compares them to results of non-stage-parallel
versions of the solvers.}
Finally, Section~\ref{sec:outlook}
summarizes our findings and points to further research directions.

\section{Fully implicit Runge--Kutta methods and stage-parallel solution approaches}\label{sec:methods}

{In the following, we summarize key aspects of the implicit Runge--Kutta methods and
stage-parallel solution procedures analyzed in this publication.
%
We consider the subclass of implicit Runge-Kutta methods referred to as Radau IIA methods. For $Q$ stages, the method order is given by $2Q-1$. 

We restrict ourselves to the case of a linear system of equations of the form
\[
M
\frac{\partial \textbf{u}(t)}{\partial t} + K \textbf{u}(t) = \textbf{g}(t),
\]
where $M$ denotes the mass matrix and $K$ the stiffness matrix.
To obtain the values of $\textbf{k}_q$ at the stages $q=1,\ldots,Q$, we need to solve the system
\[
\begin{bmatrix}
M +  \tau a_{11} K & \tau a_{12}K    &  \hdots        &    \tau a_{1Q}K     \\
 \tau a_{21}K & M + \tau a_{22}K &   \hdots       &     \tau  a_{2Q}K    \\
 \vdots & \vdots & \ddots &     \vdots      \\
\tau a_{Q1}K &\tau a_{Q2}K & \hdots  & M + \tau a_{QQ}K
\end{bmatrix} \begin{bmatrix}
\textbf{k}_1 \\ \textbf{k}_2 \\ \vdots \\ \textbf{k}_Q
\end{bmatrix} = \begin{bmatrix}
K \textbf{u}_0 + \textbf{g}(t_0+c_1 \tau) \\
K \textbf{u}_0 + \textbf{g}(t_0+c_2 \tau) \\
\vdots \\
K \textbf{u}_0 + \textbf{g}(t_0+c_Q \tau)
\end{bmatrix},
\]
which is expressed using Kronecker products as
\[
( \mathbb{I}_Q  \otimes M + \tau A_Q \otimes K) \textbf{k} = \overline{\textbf{g}} + (\mathbb{I}_Q \otimes K)(\textbf{e}_Q \otimes \textbf{u}_0 ).
\]
Next, we multiply both sides by $ (A_Q^{-1} \otimes \mathbb{I}_n) $ and use the relation $(A \otimes B)(C \otimes D)=(AC)\otimes(BD)$ to obtain
\begin{equation}\label{eq:irk:A}
\underbrace{( A_Q^{-1} \otimes M + \tau \mathbb{I}_Q \otimes K)}_{\mathcal{A}} \textbf{k} =  (A_Q^{-1} \otimes \mathbb{I}_n) \overline{\textbf{g}} + (A_Q^{-1} \otimes K)(\textbf{e}_Q \otimes \textbf{u}_0 ),
\end{equation}
which is the form we utilize in this study.
Following Butcher~\cite{butcher1976implementation}, one can construct the 
spectral decomposition of $A_Q^{-1} = S \Lambda S^{-1}$ and use this to transform the matrix:
\begin{align*}
\mathcal{A} &=
( A_Q^{-1} \otimes M + \tau \mathbb{I}_Q \otimes K )
=
(S \otimes \mathbb{I}_n) (\Lambda \otimes M + \tau \mathbb{I}_Q \otimes K)(S^{-1} \otimes \mathbb{I}_n).
\end{align*}
The inverse of the matrix and the solution
of the stages is explicitly given as:
\begin{align}\label{eq:irk:Ainv}
\textbf{k}=
\underbrace{{(S \otimes \mathbb{I}_n)} {(\Lambda \otimes M + \tau \mathbb{I}_Q \otimes K)}^{-1} {(S^{-1} \otimes \mathbb{I}_n)}}_{\mathcal{A}^{-1}}
((A_Q^{-1} \otimes \mathbb{I}_n) \overline{\textbf{g}} + (A_Q^{-1} \otimes K)(\textbf{e}_Q \otimes \textbf{u}_0 )).
\end{align}
In the context of Radau IIA methods, 
$\Lambda$ is diagonal and contains $\lfloor Q/2 \rfloor$ complex-conjugate eigenvalue pairs as well as one real
eigenvalue in the case of odd $Q$, and the matrix $S$ contains the eigenvectors. Hence, ${(\Lambda \otimes M + \tau \mathbb{I}_Q \otimes K)}$ is block-diagonal and
its inverse is given by the inverse of each block ${(\lambda_i M + \tau K)}$,
which can be computed independently. In practice, one solves blocks corresponding
to complex-conjugate eigenvalue pairs together, 
necessitating the solution of $\lfloor Q/2 \rfloor$ complex blocks via complex arithmetic or the transformation into two-by-two
real blocks~\cite{southworth2022fast} and the solution of one real block in the case of odd $Q$.

In the literature, there are additional ways to factorize $A_Q$/$A_Q^{-1}$ and to obtain
a real block system more directly: The real Schur
complement~\cite{southworth2022fast2, wanner1996solving} leads to block triangular matrices $\hat{S}$, $\hat{S}^{-1}$ and block diagonal matrix $\hat{\Lambda}$ with two-by-two blocks of the form $\begin{bmatrix} \Re(\lambda_i) & \alpha \\ -\Im(\lambda_i)^2 / \alpha & \Re(\lambda_i) \end{bmatrix}$, for an arbitrary constant $\alpha$.

Alternatively to factorizing $\mathcal{A}$ directly, one can also solve the system~\eqref{eq:irk:A} iteratively with a Krylov solver, like GMRES,
with the help of a preconditioner.
We note that for iterative solvers with suitable preconditioners, only the action of the matrix
$ (A_Q^{-1} \otimes M + \tau \mathbb{I}_Q \otimes K) $ on a vector has to be implemented, rather than the matrix itself.

Based on an observation in Axelsson~\cite{Axelsson1}, namely, that the matrix $A_Q^{-1}$ has a dominating lower triangular part, Axelsson and Neytcheva~\cite{axelsson2020numerical} proposed to decompose $A_Q^{-1} = LU$.
Here, matrix $U$ has a unit diagonal, implying that all eigenvalues of $L^{-1}A_Q^{-1}$ are equal to one, which makes $L$ suitable for constructing a preconditioner for \eqref{eq:irk:A}. In order to obtain a preconditioner allowing for stage parallelism
and real arithmetic, the spectral decomposition $L=\tilde{S} \tilde{\Lambda} \tilde{S}^{-1}$ is employed to construct the preconditioner. Its inverse is
given as:
\begin{align}\label{eq:irk:P}
P^{-1} &=
{(\tilde{S} \otimes \mathbb{I}_n)} {(\tilde{\Lambda} \otimes M + \tau \mathbb{I}_Q \otimes K)}^{-1} {(\tilde{S}^{-1} \otimes \mathbb{I}_n)}.
\end{align}
Just as before, the term $(\tilde{\Lambda} \otimes M + \tau \mathbb{I}_Q \otimes K)$ is block-diagonal and
its inverse is given by the inverse of each block ${(\tilde{\lambda}_i M + \tau K)}$. The spectral decomposition of $L$ is always real, which is the main motivation to use $L$ instead of $A_q^{-1}$. Hence, $Q$ real blocks can be solved independently, in
contrast to the complex case with only $\lceil Q/2 \rceil$
independent blocks.
An analysis of the eigenvalues of the preconditioned system is provided in~\cite{AxelDravNeyt}. In the following, we drop the tilde over $\tilde{\Lambda}$ and $\tilde{S}$,
as the meaning of these symbols is clear in the context they are used.


In \eqref{eq:irk:A}, \eqref{eq:irk:Ainv} and \eqref{eq:irk:P}, one can identify multiple independent operations:
\begin{enumerate}
\item the right-hand-side function $\textbf{g}$ can be evaluated independently for each stage,
\item the matrix-vector multiplication with the mass matrix $M$ and the stiffness matrix $K$ in $(\mathbb{I}_q \otimes M ) \textbf{k}$ and $(\mathbb{I}_Q \otimes K ) \textbf{k}$ can be performed independently for each stage\footnote{Note that the
following decomposition is applicable: $A_Q^{-1} \otimes M = (A_Q^{-1} \otimes \mathbb{I}_n) (\mathbb{I}_Q \otimes M)$.}, and
\item $Q$ or $\lceil Q/2 \rceil$ blocks involving ${(\lambda_i M + \tau \otimes K)}$ can be solved independently.
\end{enumerate}
A parallel execution across the blocks/stages is a natural choice
on modern supercomputers.

Obviously, the combination of the partial results from the stage-parallel execution via multiplication with
{$(A_Q^{-1}  \otimes \mathbb{I}_n)$}, $(S^{-1}  \otimes \mathbb{I}_n)$ or $(S  \otimes \mathbb{I}_n)$ is not independent. This step corresponds to
a linear combination of the vector or, in the latter cases, to a basis change. Such an operation
might be challenging in parallel, especially on distributed memory systems if
each stage is assigned to a distinct process, and might counteract the benefits of the parallel execution of other parts of
the algorithm.

In the following discussions of parallelization and implementation aspects, we
omit the option of direct factorizing of $\mathcal{A}$, and
consider it only in Section~\ref{sec:complex}, pointing out that the proposed concepts
are applicable in this context as well.
For the realization of a stage-parallel iterative solver including the
stage-parallel preconditioner, one needs an efficient parallel implementation of
of basic tensor operations of the form from \eqref{eq:irk:A} and \eqref{eq:irk:P},
\begin{itemize}
\item ``generalized vector scaling'' with $C \in \mathbb{R}^{n\times n}$
\begin{align}
\textbf{v}=(\mathbb{I}_Q \otimes C) \textbf{u}
\quad\leftrightarrow\quad
\textbf{v}_i = C \textbf{u}_i
\end{align}
\item ``generalized matrix-vector product'' with $D \in \mathbb{R}^{Q\times Q}$
\begin{align}
\textbf{v}=(D  \otimes \mathbb{I}_n) \textbf{u}
\quad\leftrightarrow\quad
\textbf{v}_i = \sum_{1\le j \le Q}
D_{ij}\textbf{u}_j.
\end{align}
\end{itemize}
Furthermore, we analyze the benefits of the stage--parallel solver
in comparison to its sequential counterpart. For the blocks,
we use geometric multigrid~\cite{munch2022gc}
as efficient solver with state-of-the-art parallel scaling.

We call the algorithm \eqref{eq:irk:A}+\eqref{eq:irk:P} \textit{stage-parallel IRK (in the figures abbreviated as SPIRK)} when the blocks
are solved in parallel. When the blocks are solved sequentially, we simply use \textit{IRK}.

We conclude this section by providing a brief comparison of the
preconditioner~\eqref{eq:irk:P}
and the
ones considered by Pazner and Person~\cite{pazner2017stage} solving time-dependent
non-linear problems.
They analyze
block Jacobi solvers with different block sizes: one that takes into account all
coupling terms between stages (``stage-coupled'') and on that ignores them (``stage-uncoupled''). In the latter case, the blocks corresponding to the stages are independent of
each other and can be solved in a stage-parallel way.
The structure is very similar to the inner term of~\eqref{eq:irk:P}. However, we
note that the basis changes $\tilde{S}$ and $\tilde{S}^{-1}$ imply a
coupling of the stages. We show that the cost of the basis changes is small, enabling
a stage-parallel ``stage-coupled'' preconditioner with the cost
similar to that of the stage-parallel ``stage-uncoupled'' preconditioner proposed in~\cite{pazner2017stage}.
}

\section{Implementation details}\label{sec:implementation}

In the following, we discuss pure MPI implementations of the stage-parallel IRK; in some of our experiments, we use  MPI's shared-memory features. The algorithms can be easily generalized to task-based implementations, as
provided by OpenMP, and to hybrid implementations (MPI+X).
{In Sections~\ref{sec:implementation:dd}--\ref{sec:implementation:tensor}, we
describe an approach where stages are distributed and assigned to distinct processes, and, in Section~\ref{sec:implementation:batched}, we describe a way to batch operations from different stages.
Both approaches aim to increase 
the parallelism in the solver and hereby to increase  
the sizes of the subproblems that can
be processed in parallel.}

\subsection{Domain decomposition}\label{sec:implementation:dd}
\begin{figure}
    \centering
    \begin{adjustbox}{width=0.7\textwidth} 
    \includegraphics[width=0.7\textwidth]{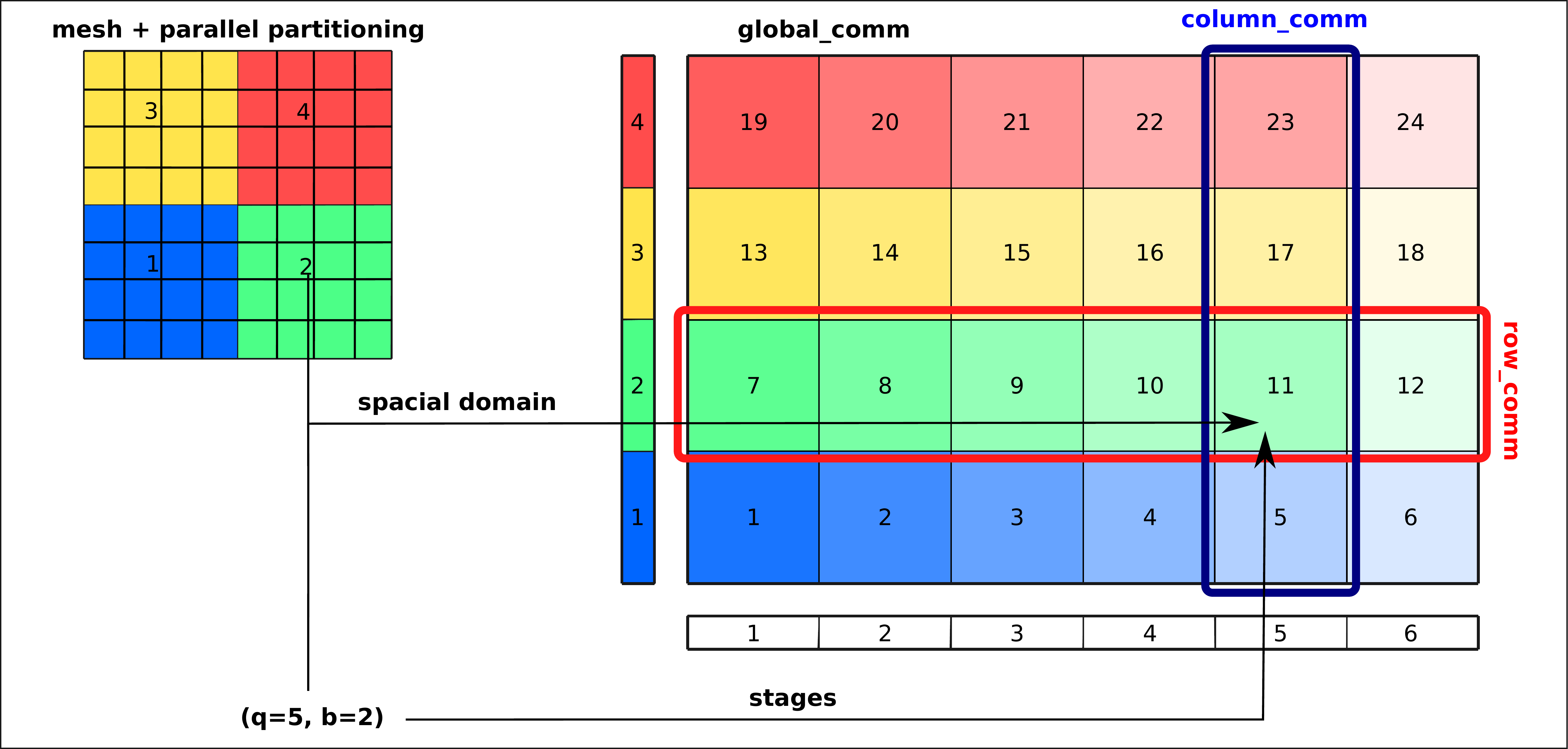}
    \end{adjustbox}
    \caption{Three MPI communicators (for a hypothetical setup with  24 processes, $Q=6$ stages, and $B=4$ partitions) used to simplify communications in stage-parallel IRK:  \texttt{global\_comm} collects processes with work and
\texttt{column\_comm}/\texttt{row\_comm} collects processes owning the same stage/partition of the computational domain.
{Furthermore, the mapping
between the stage/partition pair and the rank in the global communicator is 
indicated.} Adopted from \cite{munch2021hyper}.}\label{fig:comm}
\end{figure}
We decompose the mesh of our spatial computational domain into $B$ partitions. In the case of IRK, we assign each partition to
a (MPI) process. In the case of stage-parallel IRK, we assign each process a pair $(q,b)\in [1, Q]\times [1,B]$ consisting of a stage and a partition.
For the sake of simplicity, we enumerate the processes lexicographically, as shown in Figure~\ref{fig:comm}. However, a basic preprocessing
step based on virtual topologies allows us to use any enumeration of processes, as is needed in
later discussions.

In total, we need at least $Q\times B$ processes. In order to be able to perform operations between processes with the same stage (e.g., to solve the
inner blocks) or the same partition (e.g., for the basis change), we use
additional subcommunicators: \texttt{column\_comm} and \texttt{row\_comm}. This is a common approach in the
context of distributed matrix-matrix-multiplication implementations~\cite{van1997summa} and
finite-difference stencil computations on Cartesian meshes~\cite{hager2010introduction}. A similar approach has been used
in~\cite{munch2021hyper} to solve the 6D Vlasov--Poisson equation on the tensor product of a 3D geometric domain and a 3D velocity-space domain.

\subsection{Parallel data distribution}\label{sec:implementation:data}

A natural choice is to let each process $(q,b)$ own the locally relevant part of the
vectors associated with the stage $q$ and let it update the locally owned part of the
solution vector {during evaluation of the right-hand-side function, of the matrix-vector
product, and of the block solvers}.

For the operation $\textbf{v}=(\mathbb{I}_Q \otimes C) \textbf{u}$, setting $\textbf{v}_q=C \textbf{u}_q$ does not require communication
between stages. However, each process needs
access to the local part of the matrix $C\in \mathbb{R}^{n\times n}$ (in our case: $M$ and $K$). The need for local access to the
matrix means that certain data structures have to be duplicated $Q$ times on distributed systems. This could be addressed
by using shared memory. Since all our experiments are run without assembling any matrix (``matrix-free approach''), which is memory-efficient~\cite{kronbichler2012generic}, we defer its investigation to future work.

In contrast, the operation $\textbf{v}=(D  \otimes \mathbb{I}_n) \textbf{u}$ only needs access to a
small matrix $D \in \mathbb{R}^{Q\times Q}$, which can easily be replicated on all processes. The
major challenge for this operation is that a process needs to access the local vector entries of all
stages, e.g.,
$
\textbf{v}_2
=
D_{11} \textbf{u}_1 +
D_{21} \textbf{u}_2 +
 \dots +
D_{2Q} \textbf{u}_Q.
$
{The gathering of the needed vector entries from all stages is not feasible as this would mean that
the vectors are duplicated $Q$-times.} In the next sections, we discuss an appropriate communication
pattern to alleviate this problem and compare its performance to the one of a shared-memory approach, in which
the processes have direct read access to the needed entries of all stages.

\subsection{Distributed tensor operations}\label{sec:implementation:tensor}

To derive a memory-efficient implementation of the operation $\textbf{v}=(D  \otimes \mathbb{I}_n) \textbf{u}$, {which is needed for the linear combinations
during the setup of the right-hand-side vector and during the matrix-vector multiplication as well as for the basis changes,}
one can exploit its associativity:
\begin{scriptsize}
\begin{align}\label{eq:cannon}
\begin{bmatrix}
\textbf{v}_1 \\
\textbf{v}_2 \\
\vdots \\
\textbf{v}_Q
\end{bmatrix}
=
\underbrace{
\begin{bmatrix}
D_{11}\textbf{u}_1 \\
D_{21}\textbf{u}_1 \\
\vdots \\
D_{Q1} \textbf{u}_Q
\end{bmatrix}
+
\begin{bmatrix}
D_{12}\textbf{u}_2 \\
D_{22}\textbf{u}_2 \\
\vdots \\
D_{Q2}\textbf{u}_2
\end{bmatrix}
+
\dots
+
\begin{bmatrix}
D_{13}\textbf{u}_3 \\
D_{23}\textbf{u}_3 \\
\vdots \\
D_{Q3} \textbf{u}_3
\end{bmatrix}
}_{\bm v _i = \sum_j D_{ij}\textbf{u}_j}
=
\underbrace{
\begin{bmatrix}
D_{11}\textbf{u}_1 \\
D_{22}\textbf{u}_2 \\
\vdots \\
D_{QQ} \textbf{u}_Q
\end{bmatrix}
+
\begin{bmatrix}
D_{12}\textbf{u}_2 \\
D_{23}\textbf{u}_3 \\
\vdots \\
D_{Q1}\textbf{u}_1
\end{bmatrix}
+
\dots
+
\begin{bmatrix}
D_{13}\textbf{u}_3 \\
D_{24}\textbf{u}_4 \\
\vdots \\
D_{Q(Q-1)} \textbf{u}_{(Q-1)}
\end{bmatrix}
}_{\bm v_i = \sum_k D_{ij}\textbf{u}_{j}\text{\; with \;} j:=(i+k)\%Q}.
\end{align}
\end{scriptsize}
If we consider each summand as a computation step, we see that each process needs access to a different
part of the source vector $\textbf{u}$ associated with a stage. In order to obtain the right data access, the communication pattern is as follows.
After each computation step, the matrix is rotated to the left and the source vector upwards.
The algorithm is similar to Cannon's algorithm~\cite{cannon1969cellular}, which is
designed for distributed matrix-matrix multiplications.
In contrast
to Cannon's algorithm, we operate on tensors and with fully replicated $D$. In the case that
$\textbf{u}$ is distributed, we get a circular communication pattern, which
can be built around a sequence of calls to \texttt{MPI\_Sendrecv\_replace}, operating on \texttt{comm\_row}.

\subsection{Batching of operations}\label{sec:implementation:batched}

{Sections~\ref{sec:implementation:dd}--\ref{sec:implementation:tensor}
consider an approach that also employs parallelism across the stages by assigning each stage to a distinct compute unit. For fixed computational resources and problem sizes, this increases the size of the spatial subproblems and reduces the number of communication steps.
Alternatively, one could increase the local work by processing
$Q$ stages  on the same compute unit in a batched fashion. To be efficient, all operations of the form
$\textbf{v}=(\mathbb{I}_Q \otimes C) \textbf{u}$ in \eqref{eq:irk:A} and of
the solvers of the steps in \eqref{eq:irk:P} need to support this matrix-times-multivector mode of processing.

For the geometric multigrid solver used for the experiments in Sections~\ref{sec:performance} and \ref{sec:complex}, all ingredients in terms of smoother, prolongator/restrictor,
and coarse-grid solver need to support batching.
For simple smoothers and coarse-grid
solvers, like Chebyshev iterations around a point-Jacobi method, this is
the case, whereas it is typically not the case for algebraic multigrid.

Batching is more efficient in terms of memory consumption, since shared data
structures, particularly the matrices $M$ and $K$, only need to be stored once. Furthermore,
also data structures pertaining to matrix-free evaluation need to be loaded once per iteration. For instance,
many matrix-free implementations load the mapping data, e.g., the
Jacobian matrix of size $\mathbb{R}^{3\times3}$ in 3D, for each quadrature point during cell loops. On
affine meshes, the mapping data is the same at each quadrature point
so that it can be re-used from fast caches~\cite{kronbichler2012generic,kronbichler2019fast}. However, no simple compression is available for deformed mesh cells, where matrix-free loops become memory-bound on modern hardware. In a
sequential execution of the stage solvers, the mapping data has to be loaded from main
memory for each stage.
Also for stage parallelism according to Sections~\ref{sec:implementation:dd}--\ref{sec:implementation:tensor}, different processes load separate arrays via a shared resource, the bus from main memory, again giving a similar cost as in the stage-sequential case. In contrast, batching allows to
load the mapping data once, since all stages are processed during a single
cell loop where data from caches is still hot. Additional work for evaluating the cell integrals
of multiple stages could be hidden behind the slow memory access. 

As an alternative, on-the-fly evaluation of high-order mappings for complicated geometries is possible~\cite{kronbichler2022enhancing}, but it comes with additional complexities. We therefore concentrate on affine meshes in this work and defer the investigation of complicated meshes in the context
of stage parallelism to future work.
}

\section{Performance modeling}\label{sec:modeling}

In the following, we derive performance models for sequential and stage-parallel IRK. Since the application of the
preconditioner $P^{-1}$ is the most expensive ingredient in both cases (see the
results in Section~\ref{sec:performance}), we
consider it in detail. The statements made for the preconditioner can be straightforwardly transferred to
other parts of the algorithm.

\begin{figure}
    \centering
    \def\svgwidth{1.0\columnwidth}
    {\footnotesize 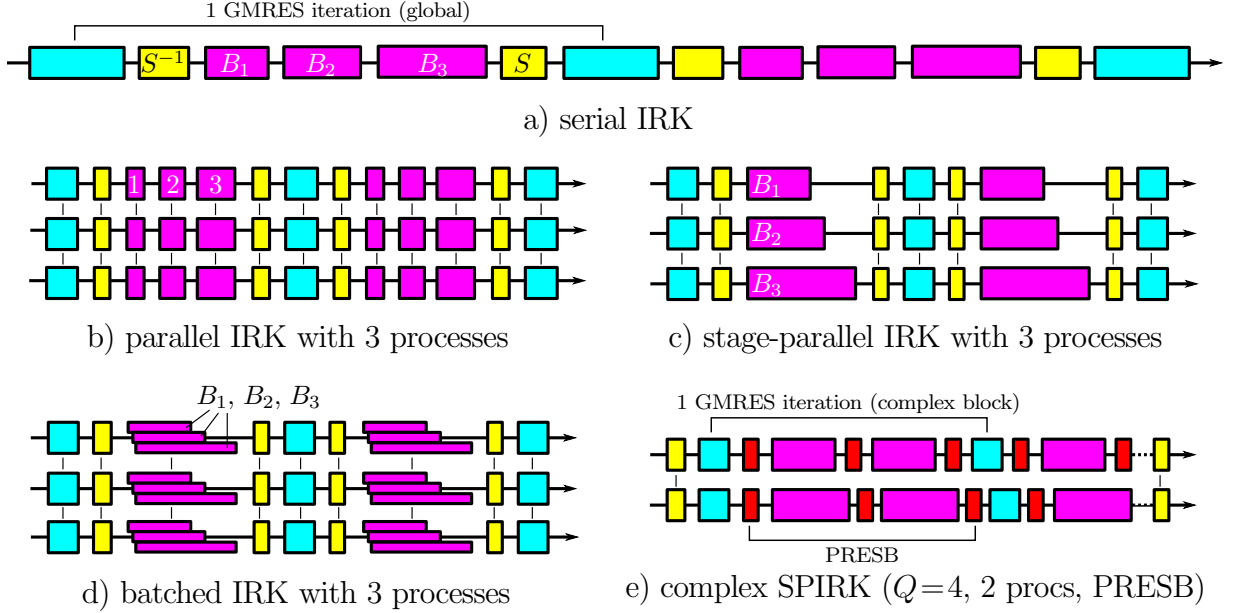}
    \caption{a-d) Visualization of a serial, a  parallel, a stage-parallel, and a batched execution of the
    application of preconditioner $P^{-1}$ of IRK with $Q=3$. The main components are the
     GMRES solver, the basis changes ($S$, $S^{-1}$), and the block solvers ($B_i$).
     {e) Visualization of a stage-parallel execution of the
    complex IRK with $Q=4$ and PRESB. Vertical lines indicate communication
    between processes.}}\label{fig:irk:trace}
\end{figure}

Figure~\ref{fig:irk:trace}a) shows a possible trace of a serial execution of IRK, and
Figure~\ref{fig:irk:trace}b) presents a parallel execution of IRK (with $3$ processes).
In both cases, the basis change $S^{-1}$, the inner block solvers $B_i$, and the basis change
$S$ are executed in sequence so that one gets for the total runtime of the
application of the preconditioner:
\begin{align}\label{eq:time:irk}
T_{\text{IRK}}(N) = 2 \cdot T_S(N) + \sum_{1 \le q \le Q} T_{B_i}(N),
\end{align}
where $T_{\text{IRK}}(N)$ denotes the total time for one application of the preconditioner. Clearly, it is a function of the number of processes $N$. In the ideal case, $\mathcal{T}_\square(N)\approx T_\square(1)/N$.
However, inherently serial parts in the code prevent perfect scaling.

In the stage-parallel case, the steps of the inner block solver are performed in parallel. Since each block may have different
properties, the solution processes might
take different amounts of time and require different numbers of inner iterations
in the case that iterative solvers are used on the blocks.
Recall that the block systems are of the form $(\lambda_i M + \tau K)$, the number of stages slightly affects the range of values of $\lambda_i$, the discretization affects the entries of $M$, and the time step $\tau$ scales the stiffness matrix $K$.

The solution
of each block is combined during the application of $S$, leading to an unavoidable synchronization.
This results in a trace, as indicated in Figure~\ref{fig:irk:trace}c), and in a total runtime that is determined by the
maximum runtime of any of the block solvers:
\begin{align}\label{eq:time:spirk}
T_{\text{SPIRK}}(N) = 2 \cdot T_S'(N) + \max_{1 \le q \le Q} \left( T_{B_q}(N/Q) \right) .
\end{align}
{Equation~\eqref{eq:time:spirk} also describes the runtime in the case of the batched approach
(see Figure~\ref{fig:irk:trace}d)).}

Based on \eqref{eq:time:irk} and \eqref{eq:time:spirk}, one can expect the following parallel performance behavior.
The timings are
comparable if 1) $T_S$ and $T_S'$ are of the same order, 2)
the block solvers are scaling nearly ideally, i.e., $T_{B_i}(N) \approx T_{B_i}(1)/N$, and 3) the solution
times of the block solvers are comparable, i.e., $\forall i\neq j:T_{B_i}\approx T_{B_j}$.
The parallel performance of the sequential IRK deteriorates if $T_\square(N)\gg T_\square(1)/N$, which is the case at
the scaling limit. The performance of the stage-parallel IRK method deteriorates if there is a significant difference in
the solution times of the block solvers, limiting the maximum speedup to $\sum_Q T_{B_q}/\max({T_{B_q}}) \le Q$.
One can deduce that stage-parallel IRK has advantages only if IRK is at the scaling
limit and the solution times of the blocks are comparable.

{In the discussion above, we assumed
wall-clock times of the solution of a block as given. However, we can refine the
expressions for iterative solvers:
\begin{align*}
T_{\text{IRK}}(N) = \sum_{1 \le q \le Q} N_i^{\text{IT}} \cdot \hat{T}_{B_i}(N)
\quad\text{and}\quad
T_{\text{SPIRK}}(N) = \max_{1 \le q \le Q} \left(  N_q^{\text{IT}}  \cdot \hat{T}_{B_q}(N/Q) \right),
\end{align*}
with $\hat{T}_i$ being the time of one iteration and $N_i^{\text{IT}}$ the number
of iterations. For the sake of simplicity, $T_S$ and $T_S'$ are dropped. If we assume that we are at the scaling
limit ($\lim_{N\to\infty}\hat{T}_q(N) \approx \lim_{N\to\infty}\hat{T}_q(N/Q)$), we get the expressions:
\begin{align*}
\lim_{N\to\infty}
T_{\text{IRK}}(N) \sim \sum_{1 \le q \le Q} N_i^{\text{IT}}
\quad\text{and}\quad
\lim_{N\to\infty}
T_{\text{SPIRK}}(N) \sim \max_{1 \le q \le Q} \left(  N_q^{\text{IT}}\right),
\end{align*}
indicating that it is possible to estimate bounds of maximum speedups based on the number of
iterations that can be run in parallel. This gives also a simple mean to compare the
benefits to alternative (stage-parallel) implementations, like the direct
factorization~\eqref{eq:irk:Ainv}, where one could consider all block solves accumulated
over all GMRES iterations for one time step.

}

In the numerical experiments in Section~\ref{sec:performance} and \ref{sec:complex}, we will evaluate the statements made above.

\section{Numerical experiments}\label{sec:performance}

\begin{table}[!t]
\centering

\caption{Number of cells and of degrees of freedom for different number of
refinements $L$.}\label{tab:results:numbers}

\begin{footnotesize}
\begin{tabular}{cccc}
\toprule
& & \multicolumn{2}{c}{degrees of freedom}\\
$L$ & cells & $k=1$ & $k=4$\\
\midrule
4	&4.1E+03	&4.9E+03	&2.7E+05 \\
5	&3.3E+04	&3.6E+04	&2.1E+06 \\
6	&2.6E+05	&2.7E+05	&1.7E+07 \\
7	&2.1E+06	&2.1E+06	&1.4E+08 \\
\bottomrule
\end{tabular}
\qquad
\begin{tabular}{cccc}
\toprule
& & \multicolumn{2}{c}{degrees of freedom}\\
$L$ & cells & $k=1$ & $k=4$\\
\midrule
8	&1.7E+07	&1.7E+07	&1.1E+09 \\
9	&1.3E+08	&1.4E+08	&8.6E+09 \\
10	&1.1E+09	&1.1E+09	&6.9E+10 \\
11	&8.6E+09	&8.6E+09	&5.5E+11 \\
\bottomrule
\end{tabular}
\end{footnotesize}
\end{table}

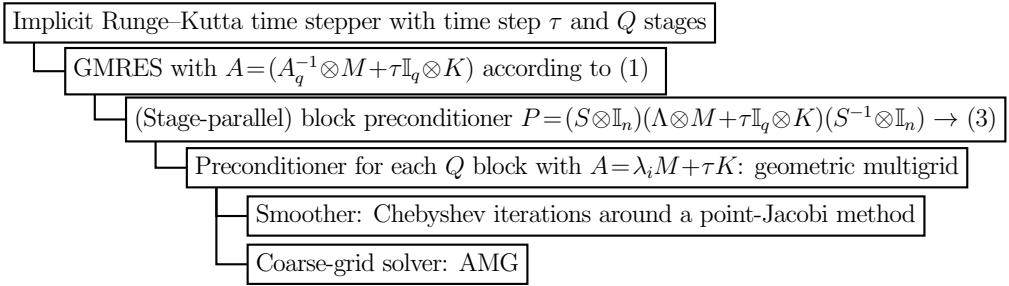
\begin{figure}[!t]

\centering
\begin{tikzpicture}[thick,scale=0.8, every node/.style={scale=0.8}]
\draw(-0.5,0.8) |- (1.0,0.0);
\draw(0.5,0) |- (1.0,-0.8);
\draw(1.5,-0.6) |- (2.0,-1.6);
\draw(2.5,-1.6) |- (3.0,-2.4);
\draw(2.5,-1.6) |- (3.0,-3.2);

\node[draw,align=left,anchor=west, fill=white] at (-1.0,0.8) {Implicit Runge--Kutta time stepper with time step $\tau$ and $Q$ stages};
\node[draw,align=left,anchor=west, fill=white] at (0,0) {GMRES with $A=( A_q^{-1} \otimes M + \tau \mathbb{I}_q \otimes K)$ according to~\Eqref{eq:irk:A} };
\node[draw,align=left,anchor=west, fill=white] at (1.0,-0.8) {(Stage-parallel) block preconditioner $P=(S \otimes \mathbb{I}_n) (\Lambda \otimes M + \tau \mathbb{I}_q \otimes K)(S^{-1} \otimes \mathbb{I}_n)$ $\rightarrow$ \Eqref{eq:irk:P}};
\node[draw,align=left,anchor=west, fill=white] at (2.0,-1.6) {Preconditioner for each $Q$ block with $A=\lambda_i M + \tau K$: geometric multigrid};
\node[draw,align=left,anchor=west, fill=white] at (3.0,-2.4) {Smoother: Chebyshev iterations around a point-Jacobi method};
\node[draw,align=left,anchor=west, fill=white] at (3.0,-3.2) {Coarse-grid solver: AMG};

\end{tikzpicture}
\caption{Diagram of the solver used to solve the heat problem in
Section~\ref{sec:performance}.}\label{fig:results:solver_diagram}
\end{figure}

In this section, we present performance results of the stage-parallel implementation
of \eqref{eq:irk:A}+\eqref{eq:irk:P}. We start with results obtained
on 16 compute nodes. In particular, we discuss the performance of a base configuration and the influence of key
parameters. We conclude the section with presenting a strong-scaling analysis.

Hereafter, we solve the 3D heat problem $\partial u/\partial t = \Delta u + f$ with the following manufactured solution:
\begin{align*}
u(x,y,z,t) = \sin(2 \pi x) \sin(2 \pi y) \sin(2 \pi z) (1 + \sin(\pi t)) \exp(-0.5 t),
\end{align*}
on a cube $\Omega = [0, 1]^3$. The source-term function $f$ and the Dirichlet boundary conditions are selected appropriately. The spatial variables are discretized with the finite element method,
for which we use the open-source library \texttt{deal.II}~\cite{dealII93, dealii2019design}. The mesh is obtained by $L$ steps of isotropic refinement of
a coarse mesh consisting of a single hexahedral cell, giving $2^{L}$ mesh cells per spatial direction or $2^{3L}$ cells in total.
We use {continuous Lagrange finite elements,
defined as the tensor products of 1D finite elements with
degree $k$}. For quadrature, we consider the consistent
Gauss--Legendre quadrature rule {with $(k+1)^3$ points}.
Table~\ref{tab:results:numbers} shows the number of cells and
the number of degrees of freedom for $k=1$ and $k=4$ for $4 \le L \le 11$.
{The time step is set to $\tau=0.1$, and we run 10 time steps.}

As outer solver of \eqref{eq:irk:A}, we apply GMRES. It is run until the $l_2$-norm of the residual has been reduced by
$10^{12}$.
As approximate inverse of each block of \eqref{eq:irk:P}, we use a single V-cycle of the geometric multigrid from \texttt{deal.II}~\cite{munch2022gc}.
As a smoother, we
apply Chebyshev iterations around a point-Jacobi method~\cite{adams2003parallel} with degree 5 and, as coarse-grid solver,
we use the algebraic multigrid solver from ML~\cite{gee2006ml}.
The corresponding solver diagram is shown in Figure~\ref{fig:results:solver_diagram}.
All operator evaluations
are performed using the matrix-free infrastructure described in~\cite{kronbichler2012generic, kronbichler2019fast} to ensure a high node-level performance, following the current trends of exascale finite-element algorithms described in~\cite{kolev2021efficient}.
Hence, we embed the IRK methods in a---with regard to communication costs---challenging context where differences are most pronounced.

All experiments
are conducted on the SuperMUC-NG supercomputer. Its compute nodes have 2 sockets (each with 24 cores of Intel Xeon
Skylake)
and the AVX-512 ISA
extension so that 8 doubles can be processed per instruction.\footnote{\url{https://doku.lrz.de/
display/PUBLIC/SuperMUC-NG}, retrieved on February 26, 2022.}
As compiler, we use \texttt{g++} (version 9.1.0) with the flags \texttt{-O3 -funroll-loops -march=skylake-avx512}.

\subsection{Moderately parallel runs}\label{sec:performance:moderatly}

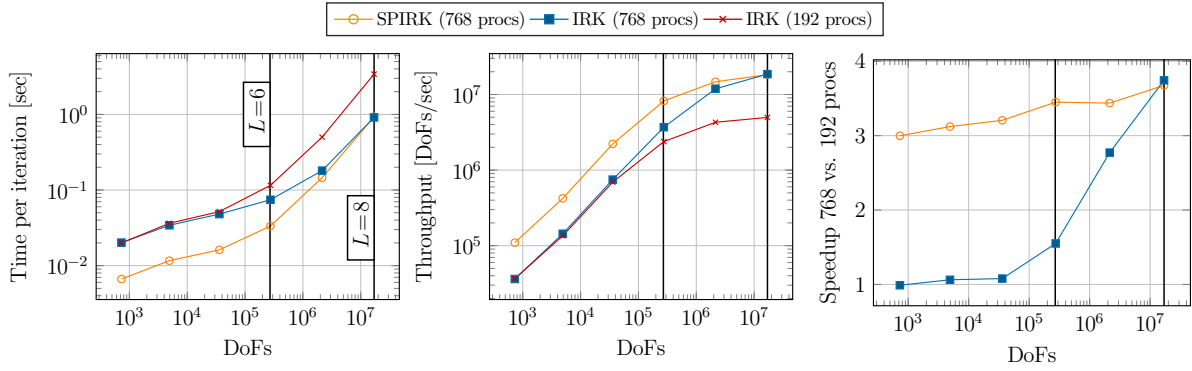
\begin{figure}

\centering

\begin{tikzpicture}[thick,scale=0.85, every node/.style={scale=0.75}]
    \begin{axis}[%
    hide axis,
    legend style={font=\small},
    xmin=10,
    xmax=1000,
    ymin=0,
    ymax=0.4,
    semithick,
    legend style={draw=white!15!black,legend cell align=left},legend columns=-1
    ]

    \addlegendimage{gnuplot@orange,every mark/.append style={fill=gnuplot@orange!80!black}, mark=o}
	    \addlegendentry{SPIRK (768 procs)};
    \addlegendimage{gnuplot@darkblue,every mark/.append style={fill=gnuplot@darkblue!80!black}, mark=square*}
	    \addlegendentry{IRK (768 procs)};
    \addlegendimage{gnuplot@red,every mark/.append style={fill=gnuplot@red!80!black}, mark=x}
	    \addlegendentry{IRK (192 procs)};

    \end{axis}
\end{tikzpicture}
\smallskip

\begin{minipage}{0.325\textwidth}
\begin{tikzpicture}[thick,scale=0.7]
    \begin{loglogaxis}[
      width=1.4\textwidth,
      height=1.2\textwidth,
      title style={font=\tiny},every axis title/.style={above left,at={(1,1)},draw=black,fill=white},
      xlabel={DoFs},
      ylabel={Time per iteration [sec]},
      legend pos={south west},
      legend cell align={left},
      cycle list name=colorGPL,
      grid,
      semithick,
      mark options={solid}
      ]

\addplot+[solid, orange,every mark/.append style={fill=gnuplot@darkblue!80!black},mark=o] table [x=dofs,y=time] {data/small-scaling/data_0.tex};
\addplot+[solid, gnuplot@darkblue,every mark/.append style={fill=gnuplot@darkblue!80!black},mark=square*] table [x=dofs,y=time] {data/small-scaling/data_1.tex};
\addplot+[solid, gnuplot@red,every mark/.append style={fill=gnuplot@red!80!black},mark=x] table [x=dofs,y=time] {data/small-scaling/data_2.tex};

\draw[draw=black,thick] (axis cs:2.7e5, 1e-3) -- node[above,sloped,inner sep=3pt,outer sep=0pt, fill=white, draw]{{  \color{black} $L=6$}}(axis cs:2.7e5, 1000);

\draw[draw=black,thick] (axis cs:1.7e7, 1e-4) -- node[above,sloped,inner sep=3pt,outer sep=0pt, fill=white, draw]{{  \color{black} $L=8$}}(axis cs:1.7e7, 10);

\end{loglogaxis}
\end{tikzpicture}
\end{minipage}
\hfill
\begin{minipage}{0.325\textwidth}
\begin{tikzpicture}[thick,scale=0.7]
    \begin{loglogaxis}[
      width=1.4\textwidth,
      height=1.2\textwidth,
      title style={font=\tiny},every axis title/.style={above left,at={(1,1)},draw=black,fill=white},
      xlabel={DoFs},
      ylabel={Throughput [DoFs/sec]},
      legend pos={south west},
      legend cell align={left},
      cycle list name=colorGPL,
      grid,
      semithick,
      mark options={solid}
      ]

\addplot+[solid, orange,every mark/.append style={fill=gnuplot@darkblue!80!black},mark=o] table [x=dofs,y=throughput] {data/small-scaling/data_0.tex};
\addplot+[solid, gnuplot@darkblue,every mark/.append style={fill=gnuplot@darkblue!80!black},mark=square*] table [x=dofs,y=throughput] {data/small-scaling/data_1.tex};
\addplot+[solid, gnuplot@red,every mark/.append style={fill=gnuplot@red!80!black},mark=x] table [x=dofs,y=throughput] {data/small-scaling/data_2.tex};

\draw[draw=black,thick] (axis cs:2.7e5, 1e4) -- node[above,sloped,inner sep=3pt,outer sep=0pt]{}(axis cs:2.7e5, 1e8);

\draw[draw=black,thick] (axis cs:1.7e7, 1e4) -- node[above,sloped,inner sep=3pt,outer sep=0pt]{}(axis cs:1.7e7, 1e8);

\end{loglogaxis}
\end{tikzpicture}
\end{minipage}
\hfill
\begin{minipage}{0.325\textwidth}
\begin{tikzpicture}[thick,scale=0.7]
    \begin{semilogxaxis}[
      width=1.45\textwidth,
      height=1.2\textwidth,
      title style={font=\tiny},every axis title/.style={above left,at={(1,1)},draw=black,fill=white},
      xlabel={DoFs},
      ylabel={Speedup 768 vs. 192 procs},
      legend pos={south west},
      legend cell align={left},
      cycle list name=colorGPL,
      grid,
      semithick,
      mark options={solid}
      ]

\addplot+[solid, orange,every mark/.append style={fill=gnuplot@darkblue!80!black},mark=o] table [x=dofs,y=speedup] {data/small-scaling/data_0.tex};
\addplot+[solid, gnuplot@darkblue,every mark/.append style={fill=gnuplot@darkblue!80!black},mark=square*] table [x=dofs,y=speedup] {data/small-scaling/data_1.tex};

\draw[draw=black,thick] (axis cs:2.7e5, 0) -- node[above,sloped,inner sep=3pt,outer sep=0pt]{}(axis cs:2.7e5, 5);

\draw[draw=black,thick] (axis cs:1.7e7, 0) -- node[above,sloped,inner sep=3pt,outer sep=0pt]{}(axis cs:1.7e7, 55);

    \end{semilogxaxis}
\end{tikzpicture}
\end{minipage}

\caption{Comparison of stage-parallel IRK with 768 processes, of IRK with 768 processes,
and of IRK with 192(=768/4) processes for $Q=4$ and $k=1$: time and throughput
per time step as well as speedup.}\label{fig:results:small:times}

\end{figure}

\begin{figure}[!t]

\centering

\begin{tikzpicture}[scale=0.44, every node/.style={scale=0.7}]
\pie[sum=auto, hide number, text = legend, explode = {0, 0, 0, 0.2, 0.2, 0.2, 0.2},radius=2.774640]
{
1.405845/rest (1.41s),
0.789860/vmult (0.79s),
0.052095/basis change (0.05s),
1.667900/block 0 (1.67s),
1.702900/block 1 (1.70s),
1.689400/block 2 (1.69s),
1.749600/block 3 (1.75s)
}
\pie[pos = {14,0}, sum=auto, hide number, text = legend, explode = {0, 0, 0, 0.2, 0.2, 0.2, 0.2},radius=3.000000]
{
2.288500/rest (2.29s),
0.960650/vmult (0.96s),
0.445950/basis change (0.45s),
6.091400/block (6.09s)
}
\end{tikzpicture}

a) $L=8$, $k=1$, $Q=4$

\begin{tikzpicture}[scale=0.44, every node/.style={scale=0.7}]
\pie[sum=auto, hide number, text = legend, explode = {0, 0, 0, 0.2, 0.2, 0.2, 0.2},radius=3.000000]
{
0.042505/rest (0.04s),
0.039612/vmult (0.04s),
0.000613/basis change (0.00s),
0.197660/block 0 (0.20s),
0.192890/block 1 (0.19s),
0.194760/block 2 (0.19s),
0.197870/block 3 (0.20s)
}
\pie[pos = {14,0}, sum=auto, hide number, text = legend, explode = {0, 0, 0, 0.2, 0.2, 0.2, 0.2},radius=1.929686]
{
0.046513/rest (0.05s),
0.021234/vmult (0.02s),
0.007543/basis change (0.01s),
0.253880/block (0.25s)
}
\end{tikzpicture}

b) $L=6$, $k=1$, $Q=4$

\caption{Time spent for matrix-vector product (vmult), basis change,
block solvers, and the remaining operations (setup of right-hand-side vector of \eqref{eq:irk:A}, vector updates, etc.) for IRK (left) and stage-parallel IRK (right). The area of the circles
indicates the total time compared to the other version.}\label{fig:results:small:pie}

\end{figure}
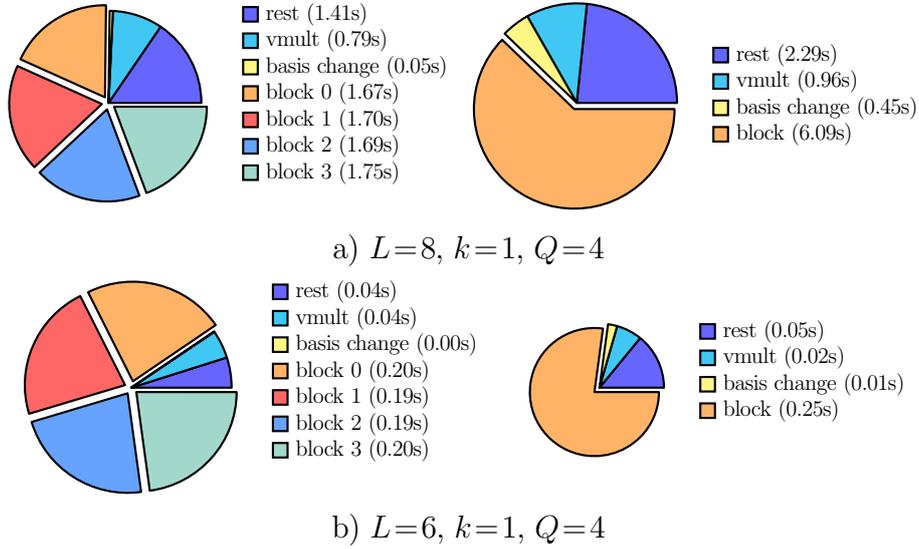

We start with the analysis of the performance of the proposed stage-parallel algorithm at small scales with 786 processes (16 compute nodes) as an example. Figure~\ref{fig:results:small:times} shows the runtime per time step and the throughput for stage-parallel IRK and IRK for $k=1$ and $Q=4$. Furthermore, we present data for
IRK with 192 (=786/4) processes, to which we compare the obtained speedup.
For large problem
sizes, IRK with four times the number of processes achieves a speedup of 3.7. 
However, its performance quickly drops for smaller sizes and times
comparable to the 192-process case are reached. In this range, the number
of processes used does not influence the times and inherently serial parts
of the code (like {latency of communication} and the coarse-grid
solver) dominate. In the
case of stage-parallel IRK, we see a different behavior: over large refinement ranges,
a speedup of $> 3$ can be reached. The maximum value is lower than in the IRK
case ($3.6 < 3.7$).

Figure~\ref{fig:results:small:pie} presents---in accordance with the traces in Figure~\ref{fig:irk:trace}---pie diagrams visualizing the time spent on different parts of the algorithms.
In particular, they show the time for basis change $S$/$S^{-1}$ and for the solution of each
block or of the whole preconditioner in the case of IRK or of stage-parallel IRK, correspondingly.
The pie diagrams are provided for two refinement configurations: $L=8$
(far away from the scaling limit) and $L=6$ (close to the scaling limit). Starting
with $L=8$, one can see that, in the case of IRK, the block preconditioners
are dominating in the total time (75\%) and the basis changes are negligible (1\%).
{In the case of stage-parallel IRK, one can
see that the time  spent on the (single) preconditioner application has decreased, but the
times for setting up the right-hand-side vector, the matrix-vector product, and
the basis changes during preconditioning have slightly increased. This is not surprising, since these are operations where communication between stages is
required and processes are implicitly synchronized by the rotation of the vectors.}
For $L=6$, the ratio of preconditioning
becomes an even more dominating part of the total solution time (90\%/77\%). However, one can see that the total time spent for preconditioning
reduces by a factor of 3.1 in the case of stage-parallel IRK compared to IRK.

\begin{figure}[!t]

\centering
\[\arraycolsep=0pt\def\arraystretch{1.5}
\begin{array}{cc}
%
%
%
\begin{minipage}{0.35\textwidth}
  \centering
  \begin{tikzpicture}[thick,scale=0.85, every node/.style={scale=0.85}]
    \begin{axis}[%
    hide axis,
    legend style={font=\footnotesize},
    xmin=10,
    xmax=1000,
    ymin=0,
    ymax=0.4,
    semithick,
    legend style={draw=white!15!black,legend cell align=left},legend columns=2
    ]

    \addlegendimage{gnuplot@darkblue,every mark/.append style={fill=gnuplot@darkblue!80!black}, mark=square*}
    \addlegendentry{$k=1$};
    \addlegendimage{gnuplot@orange, mark=diamond}
    \addlegendentry{$k=2$};
    \addlegendimage{gnuplot@red,every mark/.append style={fill=gnuplot@red!80!black}, mark=otimes}
    \addlegendentry{$k=3$};
    \addlegendimage{gnuplot@green,every mark/.append style={fill=gnuplot@green!80!black}, mark=pentagon}
    \addlegendentry{$k=4$};

    \end{axis}
  \end{tikzpicture}
  \smallskip
\end{minipage}
  &
\begin{minipage}{0.35\textwidth}
  \centering
  \begin{tikzpicture}[thick,scale=0.85, every node/.style={scale=0.85}]
    \begin{axis}[%
    hide axis,
    legend style={font=\footnotesize},
    xmin=10,
    xmax=1000,
    ymin=0,
    ymax=0.4,
    semithick,
    legend style={draw=white!15!black,legend cell align=left},legend columns=2
    ]

    \addlegendimage{gnuplot@orange,every mark/.append style={fill=gnuplot@orange!80!black}, mark=o}
    \addlegendentry{$Q=2$};
    \addlegendimage{gnuplot@darkblue,every mark/.append style={fill=gnuplot@darkblue!80!black}, mark=square*}
    \addlegendentry{$Q=4$};
    \addlegendimage{gnuplot@red,every mark/.append style={fill=gnuplot@red!80!black}, mark=x}
    \addlegendentry{$Q=6$};
    \addlegendimage{gnuplot@green,every mark/.append style={fill=gnuplot@green!80!black}, mark=triangle}
    \addlegendentry{$Q=8$};

    \end{axis}
  \end{tikzpicture}
  \smallskip
\end{minipage}
%
%
%
  \\
\begin{minipage}{0.4\textwidth}
  \flushright
 \begin{tikzpicture}[thick,scale=0.7, every node/.style={scale=0.9}]
    \begin{loglogaxis}[
      width=1.45\textwidth,
      height=1.0\textwidth,
      title style={font=\tiny},every axis title/.style={above left,at={(1,1)},draw=black,fill=white},
      xlabel={DoFs},
      ylabel={Throughput [DoFs/s/it]},
      legend pos={south west},
      legend cell align={left},
      cycle list name=colorGPL,
      grid,
      semithick,
      mark options={solid},
      xmin=300,xmax=2e8,
      ymin=5e4,ymax=1e8
      ]

\addplot+[solid, gnuplot@darkblue,every mark/.append style={fill=gnuplot@darkblue!80!black},mark=square*] table [x=dofs,y=throughput] {data/parameters/parameters-p/data_0.tex};
\addplot+[solid, gnuplot@orange,mark=diamond] table [x=dofs,y=throughput] {data/parameters/parameters-p/data_1.tex};
\addplot+[solid, gnuplot@red,every mark/.append style={fill=gnuplot@darkblue!80!black},mark=otimes] table [x=dofs,y=throughput] {data/parameters/parameters-p/data_2.tex};
\addplot+[solid, gnuplot@green,every mark/.append style={fill=gnuplot@darkblue!80!black},mark=pentagon] table [x=dofs,y=throughput] {data/parameters/parameters-p/data_3.tex};

    \end{loglogaxis}
\end{tikzpicture}
\end{minipage}
  &
\begin{minipage}{0.4\textwidth}
  \flushright
 \begin{tikzpicture}[thick,scale=0.7, every node/.style={scale=0.9}]
    \begin{loglogaxis}[
      width=1.45\textwidth,
      height=1.0\textwidth,
      title style={font=\tiny},every axis title/.style={above left,at={(1,1)},draw=black,fill=white},
      xlabel={DoFs},
      ylabel={Throughput [DoFs/s/it]},
      legend pos={south west},
      legend cell align={left},
      cycle list name=colorGPL,
      grid,
      semithick,
      mark options={solid},
      xmin=300,xmax=2e8,
      ymin=5e4,ymax=1e8
      ]

\addplot+[solid, orange,every mark/.append style={fill=gnuplot@darkblue!80!black},mark=o] table [x=dofs,y=throughput] {data/parameters/parameters-q/data_0.tex};
\addplot+[solid, gnuplot@darkblue,every mark/.append style={fill=gnuplot@darkblue!80!black},mark=square*] table [x=dofs,y=throughput] {data/parameters/parameters-q/data_1.tex};
\addplot+[solid, gnuplot@red,every mark/.append style={fill=gnuplot@darkblue!80!black},mark=x] table [x=dofs,y=throughput] {data/parameters/parameters-q/data_2.tex};
\addplot+[solid, gnuplot@green,every mark/.append style={fill=gnuplot@darkblue!80!black},mark=triangle] table [x=dofs,y=throughput] {data/parameters/parameters-q/data_3.tex};

    \end{loglogaxis}
\end{tikzpicture}
\end{minipage}
  \\
\begin{minipage}{0.4\textwidth}
  \flushright
 \begin{tikzpicture}[thick,scale=0.7, every node/.style={scale=0.9}]
    \begin{semilogxaxis}[
      width=1.45\textwidth,
      height=1.0\textwidth,
      title style={font=\tiny},every axis title/.style={above left,at={(1,1)},draw=black,fill=white},
      xlabel={DoFs},
      ylabel={Speedup vs. IRK},
      legend pos={south west},
      legend cell align={left},
      cycle list name=colorGPL,
      grid,
      semithick,
      mark options={solid},
      xmin=300,xmax=2e8
      ]

\addplot+[solid, gnuplot@darkblue,every mark/.append style={fill=gnuplot@darkblue!80!black},mark=square*] table [x=dofs,y=speedup] {data/parameters/parameters-p/data_0.tex};
\addplot+[solid, gnuplot@orange,mark=diamond] table [x=dofs,y=speedup] {data/parameters/parameters-p/data_1.tex};
\addplot+[solid, gnuplot@red,every mark/.append style={fill=gnuplot@darkblue!80!black},mark=otimes] table [x=dofs,y=speedup] {data/parameters/parameters-p/data_2.tex};
\addplot+[solid, gnuplot@green,every mark/.append style={fill=gnuplot@darkblue!80!black},mark=pentagon] table [x=dofs,y=speedup] {data/parameters/parameters-p/data_3.tex};

    \end{semilogxaxis}
\end{tikzpicture}
\end{minipage}
  &
\begin{minipage}{0.4\textwidth}
  \flushright
\begin{tikzpicture}[thick,scale=0.7, every node/.style={scale=0.9}]
    \begin{semilogxaxis}[
      width=1.45\textwidth,
      height=1.0\textwidth,
      title style={font=\tiny},every axis title/.style={above left,at={(1,1)},draw=black,fill=white},
      xlabel={DoFs},
      ylabel={Speedup vs. IRK},
      legend pos={south west},
      legend cell align={left},
      cycle list name=colorGPL,
      grid,
      semithick,
      mark options={solid},
      ytick={1,3,5},
      xmin=300,xmax=2e8,extra y ticks={2,4}
      ]

\addplot+[solid, orange,every mark/.append style={fill=gnuplot@darkblue!80!black},mark=o] table [x=dofs,y=speedup] {data/parameters/parameters-q/data_0.tex};
\addplot+[solid, gnuplot@darkblue,every mark/.append style={fill=gnuplot@darkblue!80!black},mark=square*] table [x=dofs,y=speedup] {data/parameters/parameters-q/data_1.tex};
\addplot+[solid, gnuplot@red,every mark/.append style={fill=gnuplot@darkblue!80!black},mark=x] table [x=dofs,y=speedup] {data/parameters/parameters-q/data_2.tex};
\addplot+[solid, gnuplot@green,every mark/.append style={fill=gnuplot@darkblue!80!black},mark=triangle] table [x=dofs,y=speedup] {data/parameters/parameters-q/data_3.tex};

    \end{semilogxaxis}
\end{tikzpicture}
\end{minipage}
  \\
\begin{minipage}{0.4\textwidth}
    \centering
  \footnotesize \text{(a) Vary polynomial degree $k$, fixed $Q=4$}
\end{minipage}
&
\begin{minipage}{0.4\textwidth}
    \centering
  \footnotesize \text{(b) Vary no.~stages $Q$, fixed $k=1$}
\end{minipage}
\end{array}
\]

\caption{Influence of parameters on the performance of stage-parallel IRK and
the speedup of stage-parallel IRK as compared to IRK on 16 nodes / 768 MPI processes.}\label{fig:performance:mod:parameters}

\end{figure}
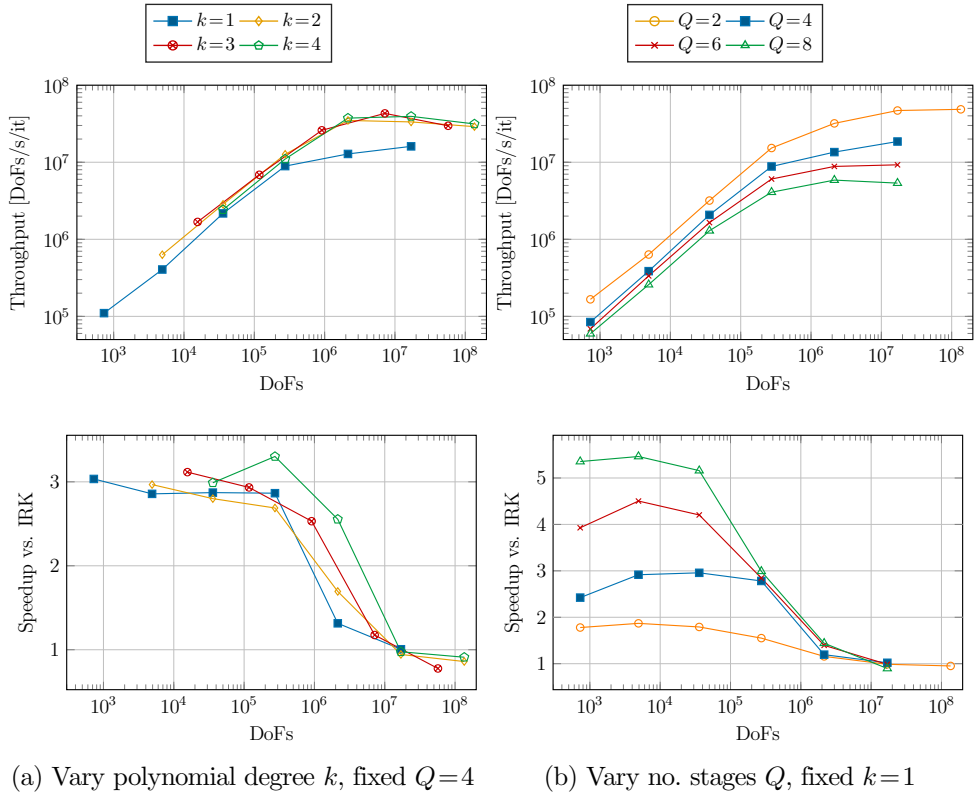

\subsubsection*{Influence of key parameters}

Figure~\ref{fig:performance:mod:parameters}
shows the influence of the parameters ``polynomial degree $k$'' and
``number of stages $Q$'' on the performance and speedup of stage-parallel IRK compared to IRK with the same
number of processes. The following
observations can be made.

Large polynomial degrees $k$ increase the throughput of
matrix-free operator evaluation.
As a consequence, the throughput of IRK also rises
with increasing polynomial degree, since its most time-consuming component---the geometric-multigrid solver---is also implemented in a matrix-free way.
Furthermore, increasing $k$ seems to also have a positive effect on
the speedup of stage-parallel IRK.

The number of stages $Q$ has the most significant influence on the throughput of IRK and stage-parallel IRK. This is not surprising, since the work per
time step increases with the number of stages. However, the observed drop in throughput
is more significant than expected from the increase in work alone. This is explained by the number of outer iterations, which increases slightly with the number of stages. The achievable speedup with
stage-parallel IRK rises with increasing number of stages. {In addition, Table~\ref{tab:complex:vs} shows the average number of GRMES iterations and accumulated numbers of $V$-cycles that are run in parallel. One can see, on the one hand, that the number of GMRES
iterations increase from 5 to 12 when going from $Q=2$ to $Q=8$ and, on the other
hand, that stage-parallel execution allows to reduce, e.g., for $Q=8$, the number
of V-cycles run  sequentially by a factor of 8 from 103 to 12, matching the measured
timings and speedups in
Figure~\ref{fig:performance:mod:parameters}.}

\begin{table}[!t]

\caption{Comparison of number of GMRES iterations (\#G) and of accumulated V-cycle applications (\#V)
  per time step for non-complex/complex IRK and stage-parallel IRK. For the stage-parallel IRK, the
  reported numbers of iterations and cycles are shown for process groups. Since
complex stage-parallel IRK runs GMRES independently on the blocks, the number of GMRES iterations and,
consequently, the number of V-cycle applications vary between blocks. In this case,
we report the minimum and maximum values. Numbers are shown for different $Q$ and $k=1$, $L=8$. The additional superscript of $\#V$ specifies the number of blocks
a GMG V-cycle considers.}\label{tab:complex:vs}

\centering

\vspace{0.2cm}

%
%
\begin{footnotesize}
\begin{tabular}{C{0.75cm}C{1.70cm}C{1.70cm}C{1.70cm}}
\hline
&\multicolumn{3}{c}{\textbf{non-complex} (Sections~\ref{sec:methods}--\ref{sec:performance})}\\ \hline
&\multicolumn{1}{c}{\textbf{IRK}}&\multicolumn{1}{c}{\textbf{SPIRK}}&\multicolumn{1}{c}{\textbf{batched}} \\ \hline
$Q$& \#G (\#V$^1$) & \#G (\#V$^1$) & \#G (\#V$^Q$) \\ \hline
2 & 5.0 (12.0) & 5.0 (6.0) & 5.0 (6.0)
\\
4 & 7.0 (32.0) & 7.0 (8.0) & 7.8 (8.8)
\\
6 & 9.9 (65.3) & 9.9 (10.9) & 9.9 (10.9)
\\
8 & 11.9 (103.1) & 11.9 (12.9) & 11.0 (12.0)
\\\hline
\end{tabular}

\vspace{0.2cm}
\begin{tabular}{C{0.75cm}C{2.00cm}C{2.0cm}C{2.70cm}C{2.5cm}}
\hline
&\multicolumn{4}{c}{\textbf{complex} (Section~\ref{sec:complex})} \\ \hline
&\multicolumn{1}{c}{\textbf{IRK-PRESB}}&\multicolumn{1}{c}{\textbf{IRK-GMG}}&\multicolumn{1}{c}{\textbf{SPIRK-PRESB}}&\multicolumn{1}{c}{\textbf{SPIRK-GMG}} \\ \hline
$Q$& \#G (\#V$^1$)& \#G (\#V$^2$)& \#G (\#V$^1$)& \#G (\#V$^2$) \\ \hline
2 & 6.0 (14.0) & 7.0 (8.0) & 6.0 (14.0) & 7.0 (8.0)
\\
4 & 11.6 (27.1) & 14.0 (16.0) & 5.6--6.0 (13.1--14.0) & 7.0--7.0 (8.0--8.0)
\\
6 & 18.2 (42.4) & 22.0 (25.0) & 6.0--6.2 (14.0--14.4) & 7.0--8.0 (8.0--9.0)
\\
8 & 23.1 (54.2) & 31.3 (35.3) & 5.0--6.1 (12.0--14.2) & 7.0--9.3 (8.0--10.3)
\\\hline
\end{tabular}
\end{footnotesize}

\end{table}

Making definite general conclusions is not straightforward as they depend on the number of processes,
the type of the coarse mesh, the refinement, the partial differential
equation, and the block solvers. However, we believe that similar
trends can be observed in different setups. We should note however, that we use the heat equation as a test problem, for which we have efficient and well-scaling block preconditioners based on multigrid methods. For other classes of problems where optimally scaling block solvers are more challenging to design, we would expect to observe differences between stage-parallel IRK and IRK also for smaller problem sizes.
{In the worst case, when the blocks are solved by direct methods,
stage parallelism might be the only way to parallelize the work. However,
there is evidence for advantages of more advanced solvers, such as the block Jacobi solver with local ILU as considered by~\cite{pazner2017stage}.}

\subsection{Virtual topology and shared memory}\label{sec:performance:vt}

\begin{figure}[!t]
    \centering
    \includegraphics[width=0.8\textwidth]{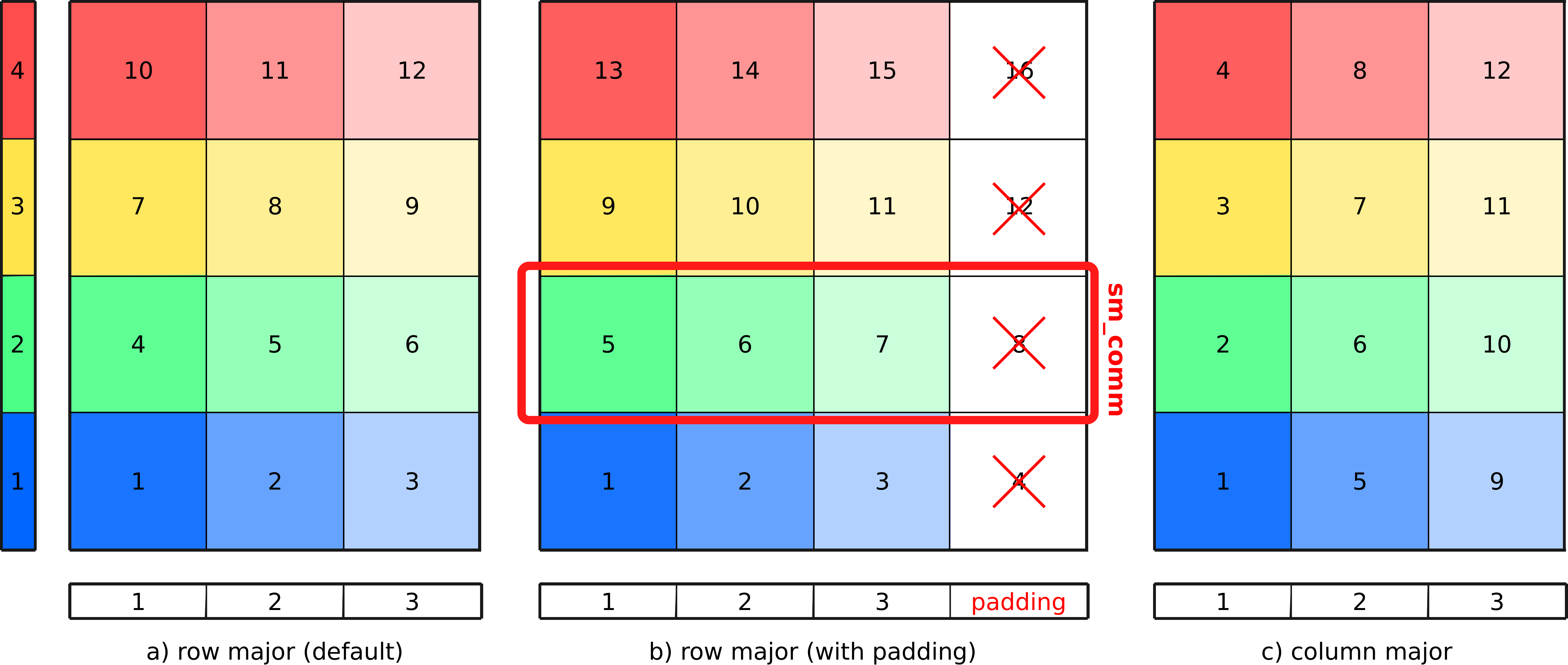}
    \caption{Different virtual topologies to increase data locality a) for the basis change $S$ and c) for the inner preconditioner. Version b) introduces padding to guarantee that all stages are on the same shared-memory domain/same compute node. 
}\label{fig:performance:vt}
\end{figure}

\begin{figure}[!t]
    \centering

\hspace{0.6cm}\begin{tikzpicture}[thick,scale=0.85, every node/.style={scale=0.85}]
    \begin{axis}[%
    hide axis,
    legend style={font=\footnotesize},
    xmin=10,
    xmax=1000,
    ymin=0,
    ymax=0.4,
    semithick,
    legend style={draw=white!15!black,legend cell align=left},legend columns=-1
    ]

    \addlegendimage{gnuplot@orange,every mark/.append style={fill=gnuplot@orange!80!black}, mark=o}
    \addlegendentry{IRK};
    \addlegendimage{gnuplot@darkblue,every mark/.append style={fill=gnuplot@darkblue!80!black}, mark=square*}
    \addlegendentry{SPIRK (row - default)};
    \addlegendimage{gnuplot@red,every mark/.append style={fill=gnuplot@red!80!black}, mark=x}
    \addlegendentry{SPIRK (row - shared memory)};
    \addlegendimage{gnuplot@green,every mark/.append style={fill=gnuplot@green!80!black}, mark=triangle}
    \addlegendentry{SPIRK (column)};

    \end{axis}
\end{tikzpicture}
\smallskip

  	\strut
  	\hfill
	 \begin{tikzpicture}[thick,scale=0.7, every node/.style={scale=0.9}]
    \begin{loglogaxis}[
      width=0.47\textwidth,
      height=0.36\textwidth,
      title style={font=\tiny},every axis title/.style={above left,at={(1,1)},draw=black,fill=white},
      title={48 procs (1 nodes)},
      xlabel={DoFs},
      ylabel={Throughput [DoFs/s/it]},
      legend pos={south west},
      legend cell align={left},
      cycle list name=colorGPL,
      grid,
      semithick,
      mark options={solid}
      ]

\addplot+[solid, orange,every mark/.append style={fill=gnuplot@darkblue!80!black},mark=o] table [x=dofs,y=throughput] {data/parameters-vt/0001/data_0.tex};
\addplot+[solid, gnuplot@darkblue,every mark/.append style={fill=gnuplot@darkblue!80!black},mark=square*] table [x=dofs,y=throughput] {data/parameters-vt/0001/data_1.tex};
\addplot+[solid, gnuplot@red,every mark/.append style={fill=gnuplot@darkblue!80!black},mark=x] table [x=dofs,y=throughput] {data/parameters-vt/0001/data_2.tex};
\addplot+[solid, gnuplot@green,every mark/.append style={fill=gnuplot@darkblue!80!black},mark=triangle] table [x=dofs,y=throughput] {data/parameters-vt/0001/data_3.tex};

    \end{loglogaxis}
\end{tikzpicture}
	\hfill
	 \begin{tikzpicture}[thick,scale=0.7, every node/.style={scale=0.9}]
    \begin{loglogaxis}[
      width=0.47\textwidth,
      height=0.36\textwidth,
      title style={font=\tiny},every axis title/.style={above left,at={(1,1)},draw=black,fill=white},
      title={768 procs (16 nodes)},
      xlabel={DoFs},
      ylabel={Throughput [DoFs/s/it]},
      legend pos={south west},
      legend cell align={left},
      cycle list name=colorGPL,
      grid,
      semithick,
      mark options={solid}
      ]

\addplot+[solid, orange,every mark/.append style={fill=gnuplot@darkblue!80!black},mark=o] table [x=dofs,y=throughput] {data/parameters-vt/0016/data_0.tex};
\addplot+[solid, gnuplot@darkblue,every mark/.append style={fill=gnuplot@darkblue!80!black},mark=square*] table [x=dofs,y=throughput] {data/parameters-vt/0016/data_1.tex};
\addplot+[solid, gnuplot@red,every mark/.append style={fill=gnuplot@darkblue!80!black},mark=x] table [x=dofs,y=throughput] {data/parameters-vt/0016/data_2.tex};
\addplot+[solid, gnuplot@green,every mark/.append style={fill=gnuplot@darkblue!80!black},mark=triangle] table [x=dofs,y=throughput] {data/parameters-vt/0016/data_3.tex};

    \end{loglogaxis}
\end{tikzpicture}
	\hfill
	 \begin{tikzpicture}[thick,scale=0.7, every node/.style={scale=0.9}]
    \begin{loglogaxis}[
      width=0.47\textwidth,
      height=0.36\textwidth,
      title style={font=\tiny},every axis title/.style={above left,at={(1,1)},draw=black,fill=white},
      title={3072 procs (64 nodes)},
      xlabel={DoFs},
      ylabel={Throughput [DoFs/s/it]},
      legend pos={south west},
      legend cell align={left},
      cycle list name=colorGPL,
      grid,
      semithick,
      mark options={solid}
      ]

\addplot+[solid, orange,every mark/.append style={fill=gnuplot@darkblue!80!black},mark=o] table [x=dofs,y=throughput] {data/parameters-vt/0064/data_0.tex};
\addplot+[solid, gnuplot@darkblue,every mark/.append style={fill=gnuplot@darkblue!80!black},mark=square*] table [x=dofs,y=throughput] {data/parameters-vt/0064/data_1.tex};
\addplot+[solid, gnuplot@red,every mark/.append style={fill=gnuplot@darkblue!80!black},mark=x] table [x=dofs,y=throughput] {data/parameters-vt/0064/data_2.tex};
\addplot+[solid, gnuplot@green,every mark/.append style={fill=gnuplot@darkblue!80!black},mark=triangle] table [x=dofs,y=throughput] {data/parameters-vt/0064/data_3.tex};

    \end{loglogaxis}
\end{tikzpicture}
	\hfill
	\strut

  	\strut
  	\hfill
	 \begin{tikzpicture}[thick,scale=0.7, every node/.style={scale=0.9}]
    \begin{semilogxaxis}[
      width=0.49\textwidth,
      height=0.32\textwidth,
      title style={font=\tiny},every axis title/.style={above left,at={(1,1)},draw=black,fill=white},
      xlabel={DoFs},
      ylabel={Speedup vs. IRK},
      legend pos={south west},
      legend cell align={left},
      cycle list name=colorGPL,
      grid,
      semithick,
      mark options={solid}
      ]

\addplot+[solid, gnuplot@darkblue,every mark/.append style={fill=gnuplot@darkblue!80!black},mark=square*] table [x=dofs,y=speedup] {data/parameters-vt/0001/data_1.tex};
\addplot+[solid, gnuplot@red,every mark/.append style={fill=gnuplot@darkblue!80!black},mark=x] table [x=dofs,y=speedup] {data/parameters-vt/0001/data_2.tex};
\addplot+[solid, gnuplot@green,every mark/.append style={fill=gnuplot@darkblue!80!black},mark=triangle] table [x=dofs,y=speedup] {data/parameters-vt/0001/data_3.tex};

    \end{semilogxaxis}
\end{tikzpicture}
  	\hfill
	 \begin{tikzpicture}[thick,scale=0.7, every node/.style={scale=0.9}]
    \begin{semilogxaxis}[
      width=0.49\textwidth,
      height=0.32\textwidth,
      title style={font=\tiny},every axis title/.style={above left,at={(1,1)},draw=black,fill=white},
      xlabel={DoFs},
      ylabel={Speedup vs. IRK},
      legend pos={south west},
      legend cell align={left},
      cycle list name=colorGPL,
      grid,
      semithick,
      mark options={solid}
      ]

\addplot+[solid, gnuplot@darkblue,every mark/.append style={fill=gnuplot@darkblue!80!black},mark=square*] table [x=dofs,y=speedup] {data/parameters-vt/0016/data_1.tex};
\addplot+[solid, gnuplot@red,every mark/.append style={fill=gnuplot@darkblue!80!black},mark=x] table [x=dofs,y=speedup] {data/parameters-vt/0016/data_2.tex};
\addplot+[solid, gnuplot@green,every mark/.append style={fill=gnuplot@darkblue!80!black},mark=triangle] table [x=dofs,y=speedup] {data/parameters-vt/0016/data_3.tex};

    \end{semilogxaxis}
\end{tikzpicture}
  	\hfill
	 \begin{tikzpicture}[thick,scale=0.7, every node/.style={scale=0.9}]
    \begin{semilogxaxis}[
      width=0.49\textwidth,
      height=0.32\textwidth,
      title style={font=\tiny},every axis title/.style={above left,at={(1,1)},draw=black,fill=white},
      xlabel={DoFs},
      ylabel={Speedup vs. IRK},
      legend pos={south west},
      legend cell align={left},
      cycle list name=colorGPL,
      grid,
      semithick,
      mark options={solid}
      ]

\addplot+[solid, gnuplot@darkblue,every mark/.append style={fill=gnuplot@darkblue!80!black},mark=square*] table [x=dofs,y=speedup] {data/parameters-vt/0064/data_1.tex};
\addplot+[solid, gnuplot@red,every mark/.append style={fill=gnuplot@darkblue!80!black},mark=x] table [x=dofs,y=speedup] {data/parameters-vt/0064/data_2.tex};
\addplot+[solid, gnuplot@green,every mark/.append style={fill=gnuplot@darkblue!80!black},mark=triangle] table [x=dofs,y=speedup] {data/parameters-vt/0064/data_3.tex};

    \end{semilogxaxis}
\end{tikzpicture}
  	\hfill
  	\strut
	\caption{Comparison of virtual topologies for $k=1$ and $Q=4$.}\label{fig:performance:vt:data}
\end{figure}
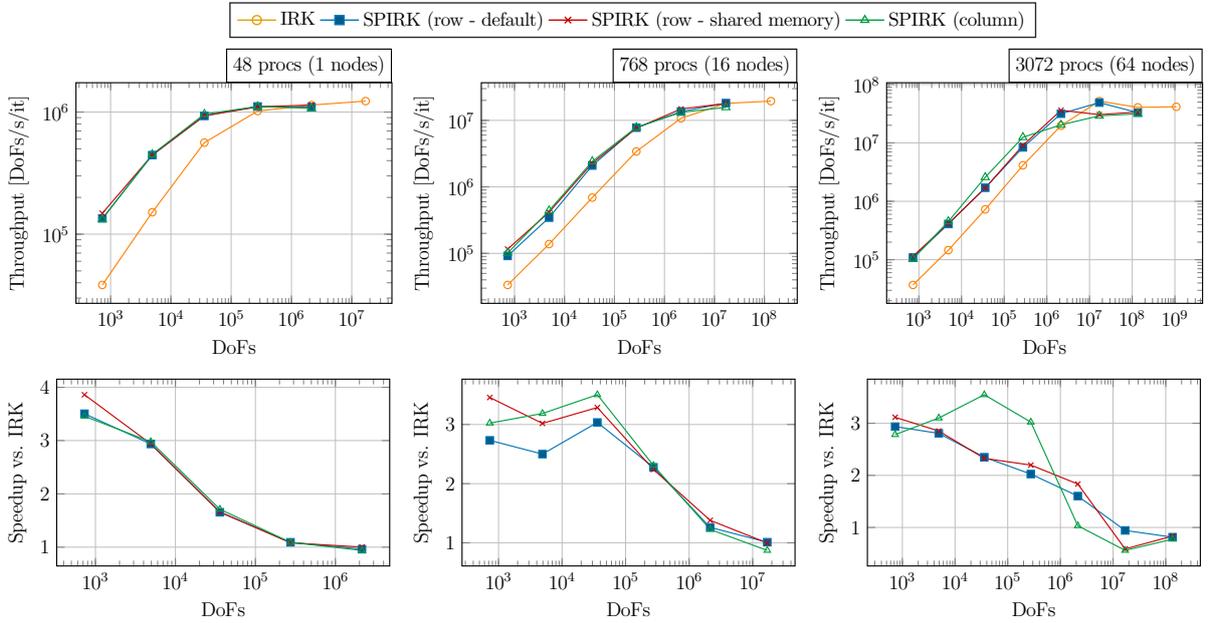

In Section~\ref{sec:implementation:dd}, we presented a row-major lexicographical enumeration of processes. This enumeration
favors the operation $\textbf{v}=(D  \otimes \mathbb{I}_n) \textbf{u}$, since the data needed are in close proximity and maybe even on the same compute node. However, it is not optimal for the inner (multigrid) solver,
for which a column-major enumeration (Figure~\ref{fig:performance:vt}c), placing nearby spatial partitions on the same node, would be favorable. In the following, we compare these two virtual
topologies.

Furthermore, we investigate the benefits of using shared-memory features of MPI
during $\textbf{v}=(D  \otimes \mathbb{I}_n) \textbf{u}$. For this, a modified row-major enumeration is beneficial.
We introduce a padding---as indicated in Figure~\ref{fig:performance:vt}b)---in such a way that all processes of
a stage are assigned to the same compute node {\color{\myred}similarly as also done
in~\cite{pazner2017stage}}. This allows us to skip the communication of complete
vectors as their entries can be accessed directly. In order to prevent race conditions, we introduce barriers across the processes
in \texttt{comm\_row}, which are equivalent to the implicit barriers in the case of parallel for loops in \texttt{OpenMP}.

Figure~\ref{fig:performance:vt:data} shows the experimental results for $k=1$ and $Q=4$ on 1, 16, and 64 compute nodes. Generally, the
different virtual topologies show similar behaviors. Using shared memory
leads to a speedup in a few cases, but overall the improvement is not significant
{in terms of performance}.
This is because the basis change is not the
bottleneck, as discussed in the previous sections. Overall,
the column-major enumeration seems to give the best results in the intermediate regime. However, it also
turns out to be the slowest virtual topology far from the
scaling limit, indicating that the increased costs of the basis change
cannot be counterbalanced by the faster block solver.

\subsection{Large-scale parallel runs}

\input{data/large-scaling/performance1}

\input{data/large-scaling/performance3}


\newcommand{\speedup}[4]{\begin{minipage}{0.242\textwidth}
\begin{tikzpicture}[thick,scale=0.76, every node/.style={scale=0.75}]
    \begin{loglogaxis}[
      width=1.6\textwidth,
      height=1.3\textwidth,
      title style={font=\tiny},every axis title/.style={above left,at={(1,1)},draw=black,fill=white},
      title={\footnotesize #2},
      xlabel={Speedup vs. IRK},
      legend pos={south west},
      legend cell align={left},
      cycle list name=colorGPL,
      grid,
      semithick,
      mark options={solid},
      ytick={#3},
      ymin={#4},
      xmin = 0.4, xmax = 9, 
      xtick = {0.5, 0.8, 1, 2, 4, 8},
      xticklabels = {0.5, \tiny{0.8}, 1, 2, 4, 8}
      ]
  
  \addplot+[only marks, red,every mark/.append style={fill=blue!80!black},mark=o] table [x=speedup,y=n_dofs] {data/large-scaling/large-scaling-2-1/speedup_#1.tex};
  \addplot+[only marks, red,every mark/.append style={fill=blue!80!black},mark=x] table [x=speedup,y=n_dofs] {data/large-scaling/large-scaling-2-4/speedup_#1.tex};

  \addplot+[only marks, blue,every mark/.append style={fill=blue!80!black},mark=x] table [x=speedup,y=n_dofs] {data/large-scaling/large-scaling-4-1/speedup_#1.tex};
  \addplot+[only marks, blue,every mark/.append style={fill=blue!80!black},mark=x] table [x=speedup,y=n_dofs] {data/large-scaling/large-scaling-4-4/speedup_#1.tex};

  \addplot+[only marks, orange,every mark/.append style={fill=blue!80!black},mark=o] table [x=speedup,y=n_dofs] {data/large-scaling/large-scaling-9-1/speedup_#1.tex};
  \addplot+[only marks, orange,every mark/.append style={fill=blue!80!black},mark=x] table [x=speedup,y=n_dofs] {data/large-scaling/large-scaling-9-4/speedup_#1.tex};


    \end{loglogaxis}
\end{tikzpicture}
\end{minipage}}

\begin{figure}
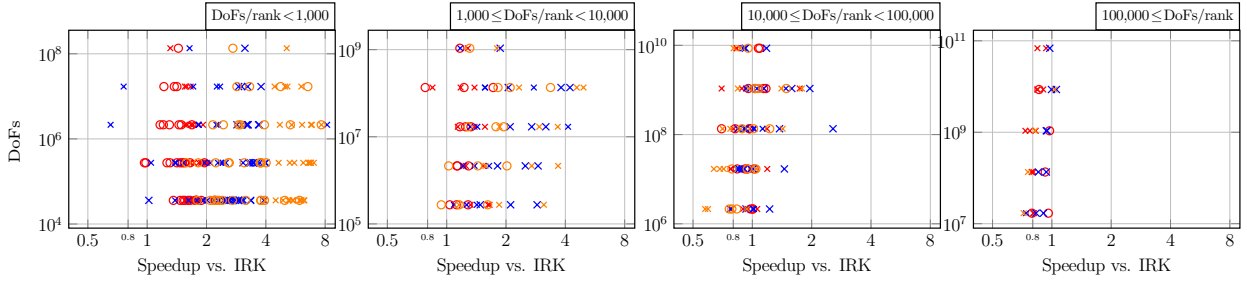


\begin{minipage}{0.015\textwidth}{\rotatebox{90}{\tiny DoFs}}\end{minipage}\,%
\speedup{0}{$\text{DoFs/rank}<1,000$}{1e4,1e6,1e8}{8e3}\,%
\speedup{1}{$1,000 \le \text{DoFs/rank}<10,000$}{1e5,1e7,1e9}{8e4}\,%
\speedup{2}{$10,000 \le \text{DoFs/rank}<100,000$}{1e6,1e8,1e10}{8e5}\,%
\speedup{3}{$100,000 \le \text{DoFs/rank}$}{1e7,1e9,1e11}{8e6}

\caption{Strong scaling: speedup categorized according to number of DoFs per process  (results from Figure~\ref{fig:performance:strongscaling}). Circles/x indicate $k=1$/$k=4$ and the colors red/blue/orange indicate $Q=2$/$Q=4$/$Q=9$.}\label{fig:resuts:speedup}
\end{figure}

Figure~\ref{fig:performance:strongscaling}
shows the results of scaling experiments starting with 1 compute node (48 processes) up to 3,072 nodes (147,456 processes) for different polynomial degrees ($k=1$/$k=4$) and numbers of stages ($Q=2$/$Q=4$/$Q=9$).
One can clearly see that stage-parallel IRK reaches lower times to solution at the
scaling limit.
Figure~\ref{fig:performance:strongscaling_normalized_all}
and \ref{fig:performance:strongscaling_normalized_stages} give more insights,
by providing normalized plots (throughput of one time step/per stage) of the same results  of all
considered values of $Q$ in a single diagram. Far from the
scaling limit (right top corner of the plots), IRK tends to be more efficient.
Stage-parallel IRK, on the contrary, reaches lower times per time step at the
scaling limit (left bottom corner of the plots) at the
cost of lower efficiencies.
Furthermore, the diagrams allow to compare the effect
of the value of $Q$, as is similarly done in Subsection~\ref{sec:performance:moderatly} for moderate number of processes.
Here, one can again see that the costs are increasing with the number
of stages, particularly also due to the increasing number of outer iterations,
which influences the scaling limits as well. However, we 
recall that a high number of stages allows to use bigger time steps due to an increased accuracy so that
the additional costs might amortize.

As a summary, Figure~\ref{fig:resuts:speedup} shows the measured speedup
of stage-parallel IRK in comparison to IRK, categorized according to the number of DoFs per process. A clear speedup
is obtained for less than 10k DoFs per process, i.e., half a million DoFs per node. For larger
problem sizes per process, the picture is split. For more than 100k DoFs per
process, IRK is consistently faster ($\approx$20\%). 

\subsection{Batched execution}\label{sec:performance:batched}

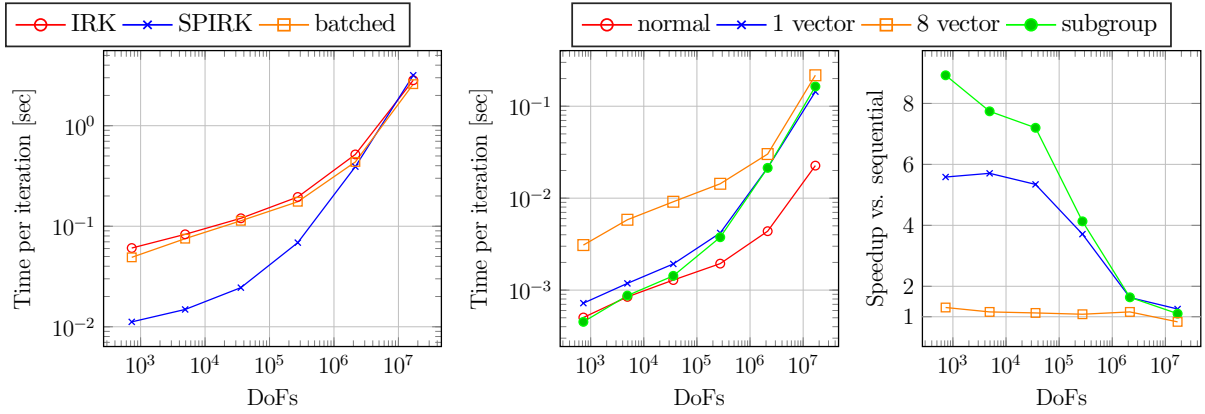
\begin{figure}

\centering

%
%
%
\begin{adjustbox}{width=\linewidth} 

\begin{tikzpicture}[thick,scale=0.85, every node/.style={scale=0.85}]
    \begin{axis}[%
    hide axis,
    legend style={font=\footnotesize},
    xmin=10,
    xmax=1000,
    ymin=0,
    ymax=0.4,
    semithick,
    legend style={draw=white!15!black,legend cell align=left},legend columns=-1
    ]
    \addlegendimage{solid, red,every mark/.append style={fill=blue!80!black},mark=o}
    \addlegendentry{IRK};
    \addlegendimage{solid, blue,every mark/.append style={fill=blue!80!black},mark=x}
    \addlegendentry{SPIRK};
    \addlegendimage{solid, orange,every mark/.append style={fill=orange!80!black},mark=square}
    \addlegendentry{batched};
    \end{axis}
\end{tikzpicture}\hspace{1.5cm}

\begin{tikzpicture}[thick,scale=0.85, every node/.style={scale=0.85}]
    \begin{axis}[%
    hide axis,
    legend style={font=\footnotesize},
    xmin=10,
    xmax=1000,
    ymin=0,
    ymax=0.4,
    semithick,
    legend style={draw=white!15!black,legend cell align=left},legend columns=-1
    ]
    \addlegendimage{solid, red,every mark/.append style={fill=blue!80!black},mark=o}
    \addlegendentry{normal};
    \addlegendimage{solid, blue,every mark/.append style={fill=blue!80!black},mark=x}
    \addlegendentry{1 vector};
    \addlegendimage{solid, orange,every mark/.append style={fill=orange!80!black},mark=square}
    \addlegendentry{8 vector};
    \addlegendimage{solid, green,every mark/.append style={fill=green!80!black},mark=*}
    \addlegendentry{subgroup};
    \end{axis}
\end{tikzpicture}\hspace{0.5cm}
 \end{adjustbox}  
\smallskip
\begin{adjustbox}{width=\linewidth} 
\strut\hfill
\begin{tikzpicture}[thick,scale=0.85, every node/.style={scale=0.9}]
  \begin{loglogaxis}[
    width=0.44\textwidth,
    height=0.4\textwidth,
    title style={font=\tiny},every axis title/.style={above left,at={(1,1)},draw=black,fill=white},
    xlabel={DoFs},
    ylabel={Time per iteration [sec]},
    legend pos={south west},
    legend cell align={left},
    cycle list name=colorGPL,
    grid,
    semithick,
    mark options={solid}
    ]

\addplot+[solid, red,every mark/.append style={fill=blue!80!black},mark=o] table [x=dofs,y=t_0] {data/parameters-all/parameters-q/data_3.tex};
\addplot+[solid, blue,every mark/.append style={fill=blue!80!black},mark=x] table [x=dofs,y=t_1] {data/parameters-all/parameters-q/data_3.tex};
\addplot+[solid, orange,every mark/.append style={fill=orange!80!black},mark=square] table [x=dofs,y=t_2] {data/parameters-all/parameters-q/data_3.tex};

  \end{loglogaxis}
\end{tikzpicture}
\hfill
 \begin{tikzpicture}[thick,scale=0.85, every node/.style={scale=0.9}]
    \begin{loglogaxis}[
      width=0.38\textwidth,
      height=0.4\textwidth,
      title style={font=\tiny},every axis title/.style={above left,at={(1,1)},draw=black,fill=white},
      xlabel={DoFs},
      ylabel={Time per iteration [sec]},
      legend pos={south west},
      legend cell align={left},
      cycle list name=colorGPL,
      grid,
      semithick,
      mark options={solid}
      ]

\addplot+[solid, red,every mark/.append style={fill=blue!80!black},mark=o] table [x=dofs,y=t] {data/gmg-batched/data_0.tex};
\addplot+[solid, blue,every mark/.append style={fill=blue!80!black},mark=x] table [x=dofs,y=t] {data/gmg-batched/data_1.tex};
\addplot+[solid, green,every mark/.append style={fill=green!80!black},mark=*] table [x=dofs,y=t] {data/gmg-batched/data_2.tex};
\addplot+[solid, orange,every mark/.append style={fill=orange!80!black},mark=square] table [x=dofs,y=t] {data/gmg-batched/data_3.tex};

    \end{loglogaxis}
\end{tikzpicture}
\hfill
 \begin{tikzpicture}[thick,scale=0.85, every node/.style={scale=0.9}]
    \begin{semilogxaxis}[
      width=0.38\textwidth,
      height=0.4\textwidth,
      title style={font=\tiny},every axis title/.style={above left,at={(1,1)},draw=black,fill=white},
      xlabel={DoFs},
      ylabel={Speedup vs. sequential},
      legend pos={south west},
      legend cell align={left},
      cycle list name=colorGPL,
      grid,
      semithick,
      mark options={solid}, extra y ticks={1},
      ]

\addplot+[solid, blue,every mark/.append style={fill=blue!80!black},mark=x] table [x=dofs,y=speedup] {data/gmg-batched/data_1.tex};
\addplot+[solid, green,every mark/.append style={fill=green!80!black},mark=*] table [x=dofs,y=speedup] {data/gmg-batched/data_2.tex};
\addplot+[solid, orange,every mark/.append style={fill=orange!80!black},mark=square] table [x=dofs,y=speedup] {data/gmg-batched/data_3.tex};

    \end{semilogxaxis}
\end{tikzpicture}
\hfill\strut
 \end{adjustbox}  

\caption{a) Comparison of time per iteration of IRK, stage-parallel IRK, and batch IRK
for $k=1$ adn $Q=8$ with 768 processes (16 compute nodes); b-c) time of iteration and speedup of different execution modes to process $Q=8$ stages in comparison to a sequential execution of the stages.}\label{fig:result:batched}

\end{figure}

{For the batched experiments, we replace the coarse-grid solver AMG
by Chebyshev iterations around a point-Jacobi method with the same settings
as of the smoothers (cf. Section~\ref{sec:performance}). Furthermore, we do not set up the coefficients of the
Chebyshev polynomials
for each block separately, but instead set them up with the approximation of the maximum
eigenvalue of all blocks. This choice does not negatively affect the number of GMRES iterations, as indicated by Table~\ref{tab:complex:vs}.

Figure~\ref{fig:result:batched}a) shows the timings of a batched execution of IRK
on 16 nodes. The results are
somewhat disappointing compared to
our performance model and the number of iterations documented in
Table~\ref{tab:complex:vs}. The timings are similar to the serial execution and only small
speedup are obtained, e.g., for $Q=8$, 5--25\%.

The cause for this behavior is that the communication during the intergrid
transfer in \texttt{deal.II} is serialized for each stage. This is due to a design choice of allowing black-box processing of the stages with a separate vector for each
stage. By collecting all stages in a single vector, this issue could be overcome, as indicated by preliminary results obtained for GMG. Figures~\ref{fig:result:batched}b-c) show
timings and speedup of a single-vector execution, of an execution with $Q$ vectors,
and of an execution with $Q$ subgroups, in comparison to a sequential execution.

The results highlight the importance
of batching all parts of the code. If this is not
possible for all identified code paths, 
the execution gets serialized and scalability is limited. The results for GMG indicate that one can expect higher speedups
(up to 5.7 for $Q=8$) for batched execution; however, the speedup will be probably
smaller than in the case of the execution on subgroups.

We point out that using non-Cartesian meshes might shift the benefit
from grouped  execution towards a batched one, since loading Jacobian
matrices per quadrature point only once and not $Q$ times (once per process group) might be beneficial, as discussed in Section~\ref{sec:implementation:batched}.}

\section{Real vs. complex arithmetic -- a comparison}\label{sec:complex}

\begin{figure}[!t]
\centering

\begin{tikzpicture}[thick,scale=0.85, every node/.style={scale=0.85}]
    \begin{axis}[%
    hide axis,
    legend style={font=\footnotesize},
    xmin=10,
    xmax=1000,
    ymin=0,
    ymax=0.4,
    semithick,
    legend style={draw=white!15!black,legend cell align=left},legend columns=4
    ]

    \addlegendimage{gnuplot@orange,every mark/.append style={fill=gnuplot@orange!80!black}, mark=o}
    \addlegendentry{$Q=2$ (IRK$\approx$SPIRK)};
    \addlegendimage{gnuplot@darkblue,every mark/.append style={fill=gnuplot@darkblue!80!black}, mark=square*}
    \addlegendentry{$Q=4$};
    \addlegendimage{gnuplot@red,every mark/.append style={fill=gnuplot@red!80!black}, mark=x}
    \addlegendentry{$Q=6$};
    \addlegendimage{gnuplot@green,every mark/.append style={fill=gnuplot@green!80!black}, mark=triangle}
    \addlegendentry{$Q=8$};

    \end{axis}
\end{tikzpicture}
\smallskip

\begin{minipage}{0.45\textwidth}
\begin{adjustbox}{width=\linewidth}
 \begin{tikzpicture}[thick,scale=0.85, every node/.style={scale=0.9}]
    \begin{loglogaxis}[
      width=1.2\textwidth,
      height=1.0\textwidth,
      title style={font=\tiny},every axis title/.style={above left,at={(1,1)},draw=black,fill=white},
      xlabel={DoFs},
      ylabel={Throughput [DoFs/s/it]},
      legend pos={south west},
      legend cell align={left},
      cycle list name=colorGPL,
      grid,
      semithick,
      mark options={solid}
      ]

\addplot+[solid, orange,every mark/.append style={fill=gnuplot@darkblue!80!black},mark=o] table [x=dofs,y=throughput] {data/parameters-complex/parameters-q/data_0.tex};
\addplot+[solid, gnuplot@darkblue,every mark/.append style={fill=gnuplot@darkblue!80!black},mark=square*] table [x=dofs,y=throughput] {data/parameters-complex/parameters-q/data_1.tex};
\addplot+[solid, gnuplot@red,every mark/.append style={fill=gnuplot@darkblue!80!black},mark=x] table [x=dofs,y=throughput] {data/parameters-complex/parameters-q/data_2.tex};
\addplot+[solid, gnuplot@green,every mark/.append style={fill=gnuplot@darkblue!80!black},mark=triangle] table [x=dofs,y=throughput] {data/parameters-complex/parameters-q/data_3.tex};

    \end{loglogaxis}
\end{tikzpicture}
\end{adjustbox}
\end{minipage}
\begin{minipage}{0.45\textwidth}
\begin{adjustbox}{width=\linewidth}
 \begin{tikzpicture}[thick,scale=0.85, every node/.style={scale=0.9}]
    \begin{semilogxaxis}[
      width=1.2\textwidth,
      height=1.0\textwidth,
      title style={font=\tiny},every axis title/.style={above left,at={(1,1)},draw=black,fill=white},
      xlabel={DoFs},
      ylabel={Speedup},
      legend pos={south west},
      legend cell align={left},
      cycle list name=colorGPL,
      grid,
      semithick,
      mark options={solid}
      ]

\addplot+[solid, orange,every mark/.append style={fill=gnuplot@darkblue!80!black},mark=o] table [x=dofs,y=speedup] {data/parameters-complex/parameters-q/data_0.tex};
\addplot+[solid, gnuplot@darkblue,every mark/.append style={fill=gnuplot@darkblue!80!black},mark=square*] table [x=dofs,y=speedup] {data/parameters-complex/parameters-q/data_1.tex};
\addplot+[solid, gnuplot@red,every mark/.append style={fill=gnuplot@darkblue!80!black},mark=x] table [x=dofs,y=speedup] {data/parameters-complex/parameters-q/data_2.tex};
\addplot+[solid, gnuplot@green,every mark/.append style={fill=gnuplot@darkblue!80!black},mark=triangle] table [x=dofs,y=speedup] {data/parameters-complex/parameters-q/data_3.tex};

    \end{semilogxaxis}
\end{tikzpicture}
\end{adjustbox}
\end{minipage}

\caption{Throughput and speedup of complex stage-parallel IRK vs. complex IRK
with PRESB for $k=1$ with 768
processes (16 compute nodes).}\label{fig:complex:throughputspeedup}

\end{figure}
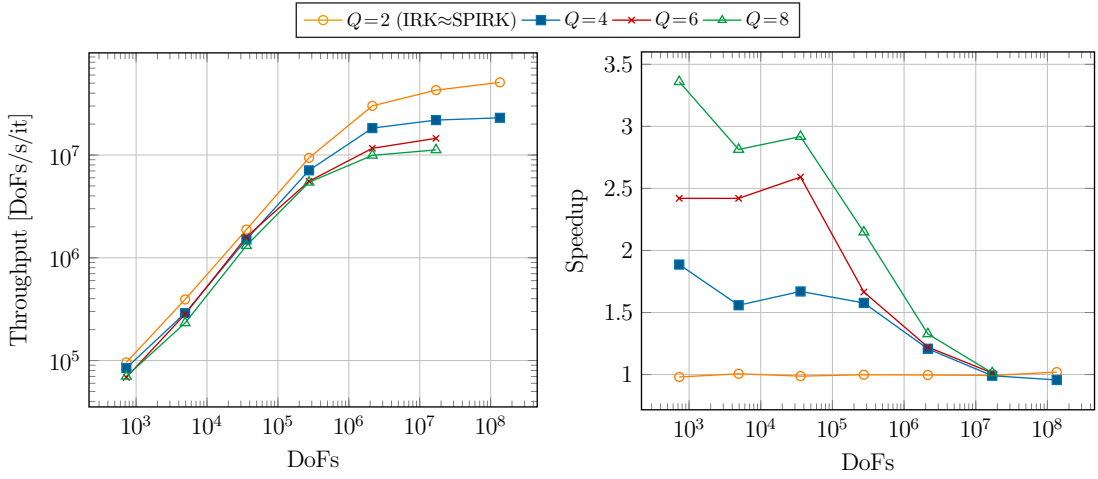

{\color{\myred}For the IRK algorithm described in Section~\ref{sec:methods}, we have chosen $L$ as basis of the preconditioner,
since its eigenvalues and eigenvectors are real and the resulting block
system to be solved is real as well. However, the proposed concepts regarding stage-parallel
implementation---particularly the data distribution and the communication
patterns---are also applicable to the complex case, which arises when $A_Q^{-1}$
is diagonalized directly. The advantage of factorizing $A_Q^{-1}$ directly is
that no global GMRES iterations---consisting of all $Q$ stages---are needed and one can solve each block individually
after the basis change. However, the disadvantage is that each block in \eqref{eq:irk:Ainv}
is complex:
\begin{align}\label{eq:presb:a}
\lambda_i M + \tau K
=
(\Re(\lambda_i) + i \cdot \Im(\lambda_i)) M + \tau K
=
\underbrace{(\Re(\lambda_i) M + \tau K)}_{K'_i} + i \cdot \underbrace{\Im(\lambda_i) M}_{M'_i},
\end{align}
with $\lambda_i = \Re(\lambda_i) + i \cdot \Im(\lambda_i)$.
Written as a two-by-two block-matrix system, $(K_i' + i M_i') \bm u_i = \bm v_i$ becomes:
\begin{align}\label{eq:complex:twobytwo}
\begin{bmatrix}
K'_i & -M'_i \\
M'_i & K_i
\end{bmatrix}
\begin{bmatrix}
\Re(\bm u_i) \\
\Im(\bm u_i)
\end{bmatrix}
=
\begin{bmatrix}
\Re(\bm v_i) \\
\Im(\bm v_i)
\end{bmatrix}.
\end{align}
Since the structure of the resulting system is similar to that of a real Schur
complement and, therefore, the algorithms to solve the system are equivalent, we will
not detail this approach.

For systems of the type~\eqref{eq:complex:twobytwo}, PRESB is an efficient preconditioner~\cite{axelsson2021robust}:
\begin{align}\label{eq:presb:p}
P_i =
\begin{bmatrix}
K'_i & -M'_i \\
M'_i & K'_i + 2 M'_i
\end{bmatrix}
=
\begin{bmatrix}
I & -I \\
0 & I
\end{bmatrix}
\begin{bmatrix}
K'_i + M'_i & 0 \\
M'_i & K'_i + M'_i
\end{bmatrix}
\begin{bmatrix}
I & I \\
0 & I
\end{bmatrix}
\end{align}
For approximating the inverse of $(K'_i+M'_i)$, we use a single V-cycle.
The corresponding solver diagram is shown in Figure~\ref{fig:presb:solver_diagram}a).
Please note that the blocks corresponding to the real and imaginary part have to
be solved in sequence (2 V-cycles), possibly limiting the scalability.

\begin{figure}[!t]

\centering

\begin{tikzpicture}[thick,scale=0.8, every node/.style={scale=0.8}]
\draw(-0.5,0.8) |- (1.0,0.0);
\draw(0.5,0) |- (0.5,-0.6);

\draw(0.5,-0.8) |- (1.0,-1.3);
\draw(1.5,-1.3) |- (2.0,-2.1);
\draw(2.5,-2.1) |- (3.0,-2.9);
\draw(2.5,-2.9) |- (3.0,-3.7);
\draw(1.5,-1.3) |- (4.0,-4.5);

\node[align=left,anchor=west] at (-0.2,-1.1) {a)~};

\node[draw,align=left,anchor=west, fill=white] at (-1.0,0.8) {Implicit Runge--Kutta time stepper with time step $\tau$ and $Q$ stages};
\node[draw,align=left,anchor=west, fill=white] at (0,0) {(Stage-parallel) GMRES for each $\lceil Q/2 \rceil$ complex block with $A=(\Re(\lambda_i) M + \tau K) + i \cdot \Im(\lambda_i) M$ $\rightarrow$ \Eqref{eq:presb:a} };

\node[draw,align=left,anchor=west, fill=white] at (1.0,-1.3) {Preconditioner: PRESB according to \Eqref{eq:presb:p}};
\node[draw,align=left,anchor=west, fill=white] at (2.0,-2.1) {Preconditioner for block 1: geometric multigrid for $( (\Re(\lambda_i) + \Im(\lambda_i)) M + \tau K )$};
\node[draw,align=left,anchor=west, fill=white] at (3.0,-2.9) {Smoother: Chebyshev iterations around a point-Jacobi method};
\node[draw,align=left,anchor=west, fill=white] at (3.0,-3.7) {Coarse-grid solver: AMG};
\node[draw,align=left,anchor=west, fill=white] at (2.0,-4.5) {Preconditioner for block 2: geometric multigrid (same as for block 1)};

\draw(0.5,-5.1) |- (1.0,-5.6);
\draw(1.5,-5.6) |- (2.0,-6.4);

\node[align=left,anchor=west] at (-0.2,-5.3) {b)};

\node[draw,align=left,anchor=west, fill=white] at (1.0,-5.6) {Preconditioner: geometric multigrid};
\node[draw,align=left,anchor=west, fill=white] at (2.0,-6.4) {Smoother/Coarse-grid solver: Chebyshev iterations around a point-Jacobi method};

\end{tikzpicture}
\caption{Diagram of the complex solvers used to solve the heat problem in
Section~\ref{sec:complex}: a) PRESB as preconditioner, b) GMG as preconditioner.}\label{fig:presb:solver_diagram}
\end{figure}
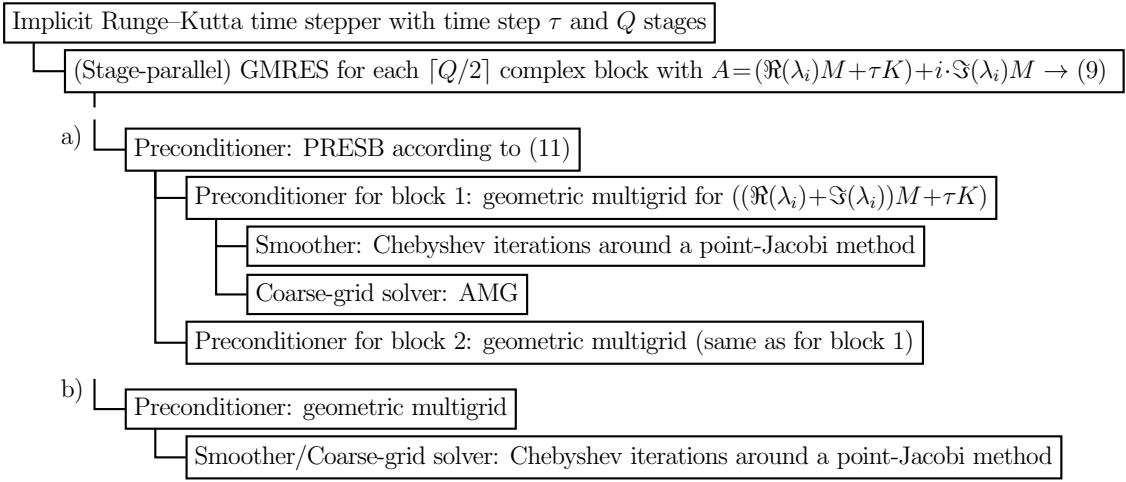

%

{Alternatively to PRESB, we also consider GMG applied directly to
\eqref{eq:complex:twobytwo}, consisting of both real and imaginary
blocks. The corresponding solver diagram is shown in Figure~\ref{fig:presb:solver_diagram}b).
The motivation for this is similar to the one of batched IRK: at the scaling
limit, running one V-cycle on a vector with twice as many DoFs might be cheaper
than running 2 V-cycles on smaller vectors in sequence.}

{The block structure of the transformed system~\eqref{eq:irk:Ainv}} has some influence on the algorithms presented in
Section~\ref{sec:implementation}.
The fact that only $\lceil Q/2 \rceil$ blocks can be solved independently
allows to parallelize only between $\lceil Q/2 \rceil$ ``stages''. Naturally, each
process group would be responsible for a stage pair. Furthermore,
operations of the form $\textbf{v}=(D  \otimes \mathbb{I}_n) \textbf{u}$
need minor adjustments, depending on the usage:
\begin{itemize}
\item application of $A_Q^{-1}$ ($\bm u, \bm v \in \mathbb{R}^{n \times Q}$, $D \in \mathbb{R}^{Q\times Q}$) with stages assigned to $\lceil Q/2 \rceil$ processes and
\item application of $S^{-1}$ ($\bm u \in \mathbb{R}^{n \times Q}$, $\bm v \in \mathbb{C}^{n \times \lceil Q/2 \rceil}$, $D \in \mathbb{C}^{\lceil Q/2 \rceil\times Q}$) and of  $S$ ($\bm u \in \mathbb{C}^{n \times \lceil Q/2 \rceil}$, $\bm v \in \mathbb{R}^{n \times Q}$, $D \in \mathbb{C}^{Q\times \lceil Q/2 \rceil}$).
\end{itemize}
In these cases, the distributed tensor algorithms from Section~\ref{sec:implementation} can be extended. They work on
blocks of two vectors and perform block operations between the rotation steps, e.g.,
\begin{align*}
\begin{bmatrix}
\textbf{v}_1 \\
\textbf{v}_2 \\  \hdashline[3pt/3pt]
\vdots
\end{bmatrix}
=
\begin{bmatrix}
D_{11}\textbf{u}_1 + D_{11}\textbf{u}_2 \\
D_{21}\textbf{u}_1 + D_{22}\textbf{u}_2 \\ \hdashline[3pt/3pt]
\vdots
\end{bmatrix}
+
\begin{bmatrix}
D_{13}\textbf{u}_3 + D_{14}\textbf{u}_4 \\
D_{23}\textbf{u}_3 + D_{24}\textbf{u}_4 \\ \hdashline[3pt/3pt]
\vdots
\end{bmatrix}
+
\dots
\end{align*}
as block extension of~\eqref{eq:cannon}.

The performance models derived in Section~\ref{sec:modeling} are
also applicable to the complex case, with more pressure on the block solvers/preconditioners. The basis change---resulting in an implicit
synchronization between all stages---has to be performed only once
per time step, since the $\lceil Q/2 \rceil$ blocks can be solved
independently. Each block solve might be more expensive, since they
might involve the solution of a two-by-two system. Comparing the complex case with
the approximate case in real arithmetic, it is obvious that we lose some possibilities for (stage)
parallelization by the reduction of the number of blocks
($Q$ $\rightarrow$ $\lceil Q/2 \rceil$) and by the sequential execution of two GMG V-cycles in the PRESB case. } {Not surprisingly, the obtained speedup
(Figure~\ref{fig:complex:throughputspeedup}) is about half of the one of the
stage-parallel preconditioner, which allows the parallel execution of all stages (see Figure~\ref{fig:performance:mod:parameters}). However, only two basis changes need to
be performed per time step, which reduces the number of possibly
expensive synchronization points.}

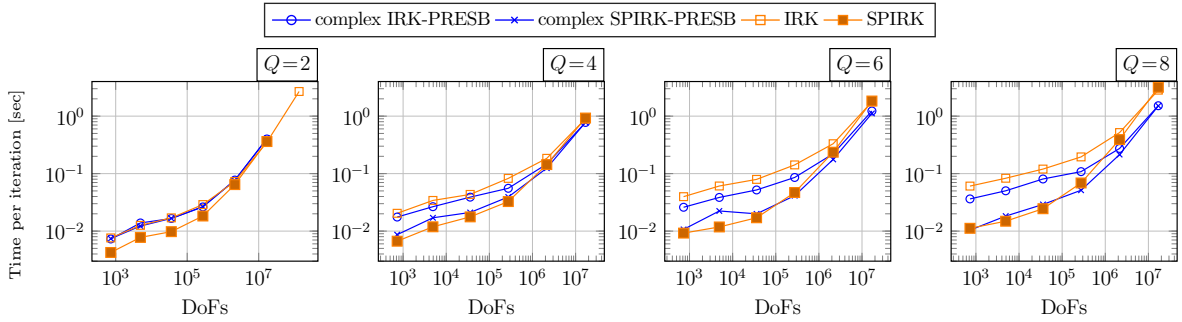
\begin{figure}[!t]
\centering

%
%
%
%
%

\begin{tikzpicture}[thick,scale=0.85, every node/.style={scale=0.85}]
    \begin{axis}[%
    hide axis,
    legend style={font=\footnotesize},
    xmin=10,
    xmax=1000,
    ymin=0,
    ymax=0.4,
    semithick,
    legend style={draw=white!15!black,legend cell align=left},legend columns=4
    ]

    \addlegendimage{solid, blue,every mark/.append style={fill=blue!80!black},mark=o}
    \addlegendentry{complex IRK-PRESB};
    \addlegendimage{solid, blue,every mark/.append style={fill=blue!80!black},mark=x}
    \addlegendentry{complex SPIRK-PRESB};
    \addlegendimage{solid, orange,every mark/.append style={fill=orange!80!black},mark=square}
    \addlegendentry{IRK};
    \addlegendimage{solid, orange,every mark/.append style={fill=orange!80!black},mark=square*}
    \addlegendentry{SPIRK};

    \end{axis}
\end{tikzpicture}
\smallskip

\begin{minipage}{0.02\textwidth}
  \tiny\rotatebox{90}{Time per iteration [sec]}
\end{minipage}
\hspace{-0.2cm}
\begin{minipage}{0.25\textwidth}
 \begin{tikzpicture}[thick,scale=0.74, every node/.style={scale=0.9}]
    \begin{loglogaxis}[
      width=1.4\textwidth,
      height=1.2\textwidth,
      title style={font=\tiny},every axis title/.style={above left,at={(1,1)},draw=black,fill=white},
      title={$Q=2$},
      xlabel={DoFs},
      legend pos={south west},
      legend cell align={left},
      cycle list name=colorGPL,
      grid,
      semithick,
      ymin={3e-3},ymax={4},
      mark options={solid}
      ]

      \addplot+[solid, blue,every mark/.append style={fill=blue!80!black},mark=o] table [x=dofs,y=t_3] {data/parameters-all/parameters-q/data_0.tex};
      \addplot+[solid, blue,every mark/.append style={fill=blue!80!black},mark=x] table [x=dofs,y=t_4] {data/parameters-all/parameters-q/data_0.tex};
      \addplot+[solid, orange,every mark/.append style={fill=orange!80!black},mark=square] table [x=dofs,y=t_0] {data/parameters-all/parameters-q/data_0.tex};
      \addplot+[solid, orange,every mark/.append style={fill=orange!80!black},mark=square*] table [x=dofs,y=t_1] {data/parameters-all/parameters-q/data_0.tex};

    \end{loglogaxis}
\end{tikzpicture}
\end{minipage}
\hspace{-0.5cm}
\begin{minipage}{0.25\textwidth}
\begin{tikzpicture}[thick,scale=0.74, every node/.style={scale=0.9}]
  \begin{loglogaxis}[
    width=1.4\textwidth,
    height=1.2\textwidth,
    title style={font=\tiny},every axis title/.style={above left,at={(1,1)},draw=black,fill=white},
    title={$Q=4$},
    xlabel={DoFs},
    legend pos={south west},
    legend cell align={left},
    cycle list name=colorGPL,
    grid,
    semithick,
    ymin={3e-3},ymax={4},
    mark options={solid}
    ]

    \addplot+[solid, blue,every mark/.append style={fill=blue!80!black},mark=o] table [x=dofs,y=t_3] {data/parameters-all/parameters-q/data_1.tex};
    \addplot+[solid, blue,every mark/.append style={fill=blue!80!black},mark=x] table [x=dofs,y=t_4] {data/parameters-all/parameters-q/data_1.tex};
    \addplot+[solid, orange,every mark/.append style={fill=orange!80!black},mark=square] table [x=dofs,y=t_0] {data/parameters-all/parameters-q/data_1.tex};
    \addplot+[solid, orange,every mark/.append style={fill=orange!80!black},mark=square*] table [x=dofs,y=t_1] {data/parameters-all/parameters-q/data_1.tex};

  \end{loglogaxis}
\end{tikzpicture}
\end{minipage}
\hspace{-0.5cm}
\begin{minipage}{0.25\textwidth}
\begin{tikzpicture}[thick,scale=0.74, every node/.style={scale=0.9}]
  \begin{loglogaxis}[
    width=1.4\textwidth,
    height=1.2\textwidth,
    title style={font=\tiny},every axis title/.style={above left,at={(1,1)},draw=black,fill=white},
    title={$Q=6$},
    xlabel={DoFs},
    legend pos={south west},
    legend cell align={left},
    cycle list name=colorGPL,
    grid,
    semithick,
    ymin={3e-3},ymax={4},
    mark options={solid}
    ]

    \addplot+[solid, blue,every mark/.append style={fill=blue!80!black},mark=o] table [x=dofs,y=t_3] {data/parameters-all/parameters-q/data_2.tex};
    \addplot+[solid, blue,every mark/.append style={fill=blue!80!black},mark=x] table [x=dofs,y=t_4] {data/parameters-all/parameters-q/data_2.tex};
    \addplot+[solid, orange,every mark/.append style={fill=orange!80!black},mark=square] table [x=dofs,y=t_0] {data/parameters-all/parameters-q/data_2.tex};
    \addplot+[solid, orange,every mark/.append style={fill=orange!80!black},mark=square*] table [x=dofs,y=t_1] {data/parameters-all/parameters-q/data_2.tex};

  \end{loglogaxis}
\end{tikzpicture}
\end{minipage}
\hspace{-0.5cm}
\begin{minipage}{0.25\textwidth}
\begin{tikzpicture}[thick,scale=0.74, every node/.style={scale=0.9}]
  \begin{loglogaxis}[
    width=1.4\textwidth,
    height=1.2\textwidth,
    title style={font=\tiny},every axis title/.style={above left,at={(1,1)},draw=black,fill=white},
    title={$Q=8$},
    xlabel={DoFs},
    legend pos={south west},
    legend cell align={left},
    cycle list name=colorGPL,
    grid,
    semithick,
    ymin={3e-3},ymax={4},
    mark options={solid}
    ]

    \addplot+[solid, blue,every mark/.append style={fill=blue!80!black},mark=o] table [x=dofs,y=t_3] {data/parameters-all/parameters-q/data_3.tex};
    \addplot+[solid, blue,every mark/.append style={fill=blue!80!black},mark=x] table [x=dofs,y=t_4] {data/parameters-all/parameters-q/data_3.tex};
    \addplot+[solid, orange,every mark/.append style={fill=orange!80!black},mark=square] table [x=dofs,y=t_0] {data/parameters-all/parameters-q/data_3.tex};
    \addplot+[solid, orange,every mark/.append style={fill=orange!80!black},mark=square*] table [x=dofs,y=t_1] {data/parameters-all/parameters-q/data_3.tex};

  \end{loglogaxis}
\end{tikzpicture}
\end{minipage}

\caption{Comparison of throughputs of non-complex/complex IRK and stage-parallel IRK for $k=1$ with 768
  processes (16 compute nodes).}\label{fig:complex:noncomplexvscomplex}
\end{figure}

\begin{figure}
\centering
\begin{tikzpicture}[thick,scale=0.85, every node/.style={scale=0.85}]
    \begin{axis}[%
    hide axis,
    legend style={font=\footnotesize},
    xmin=10,
    xmax=1000,
    ymin=0,
    ymax=0.4,
    semithick,
    legend style={draw=white!15!black,legend cell align=left},legend columns=4
    ]

    \addlegendimage{solid, blue,every mark/.append style={fill=blue!80!black},mark=o}
    \addlegendentry{complex IRK-PRESB};
    \addlegendimage{solid, blue,every mark/.append style={fill=blue!80!black},mark=x}
    \addlegendentry{complex SPIRK-PRESB};
    \addlegendimage{solid, orange,every mark/.append style={fill=orange!80!black},mark=square}
    \addlegendentry{complex IRK-GMG};
    \addlegendimage{solid, orange,every mark/.append style={fill=orange!80!black},mark=square*}
    \addlegendentry{complex SPIRK-GMG};

    \end{axis}
\end{tikzpicture}
\smallskip

\begin{minipage}{0.02\textwidth}
  \tiny\rotatebox{90}{Time per iteration [sec]}
\end{minipage}
\hspace{-0.2cm}
\begin{minipage}{0.25\textwidth}
\begin{tikzpicture}[thick,scale=0.74, every node/.style={scale=0.9}]
  \begin{loglogaxis}[
    width=1.4\textwidth,
    height=1.2\textwidth,
    title style={font=\tiny},every axis title/.style={above left,at={(1,1)},draw=black,fill=white},
    title={$Q=2$},
    xlabel={DoFs},
    legend pos={south west},
    legend cell align={left},
    cycle list name=colorGPL,
    grid,
    semithick,
    ymin={5e-3},ymax={2},
    mark options={solid}
    ]

      \addplot+[solid, blue,every mark/.append style={fill=blue!80!black},mark=o] table [x=dofs,y=t_3] {data/parameters-all/parameters-q/data_0.tex};
      \addplot+[solid, blue,every mark/.append style={fill=blue!80!black},mark=x] table [x=dofs,y=t_4] {data/parameters-all/parameters-q/data_0.tex};
      \addplot+[solid, orange,every mark/.append style={fill=orange!80!black},mark=square] table [x=dofs,y=t_5] {data/parameters-all/parameters-q/data_0.tex};
      \addplot+[solid, orange,every mark/.append style={fill=orange!80!black},mark=square*] table [x=dofs,y=t_6] {data/parameters-all/parameters-q/data_0.tex};

    \end{loglogaxis}
\end{tikzpicture}
\end{minipage}
\hspace{-0.5cm}
\begin{minipage}{0.25\textwidth}
\begin{tikzpicture}[thick,scale=0.74, every node/.style={scale=0.9}]
  \begin{loglogaxis}[
    width=1.4\textwidth,
    height=1.2\textwidth,
    title style={font=\tiny},every axis title/.style={above left,at={(1,1)},draw=black,fill=white},
    title={$Q=4$},
    xlabel={DoFs},
    legend pos={south west},
    legend cell align={left},
    cycle list name=colorGPL,
    grid,
    semithick,
    ymin={5e-3},ymax={2},
    mark options={solid}
    ]

    \addplot+[solid, blue,every mark/.append style={fill=blue!80!black},mark=o] table [x=dofs,y=t_3] {data/parameters-all/parameters-q/data_1.tex};
    \addplot+[solid, blue,every mark/.append style={fill=blue!80!black},mark=x] table [x=dofs,y=t_4] {data/parameters-all/parameters-q/data_1.tex};
    \addplot+[solid, orange,every mark/.append style={fill=orange!80!black},mark=square] table [x=dofs,y=t_5] {data/parameters-all/parameters-q/data_1.tex};
    \addplot+[solid, orange,every mark/.append style={fill=orange!80!black},mark=square*] table [x=dofs,y=t_6] {data/parameters-all/parameters-q/data_1.tex};

  \end{loglogaxis}
\end{tikzpicture}
\end{minipage}
\hspace{-0.5cm}
\begin{minipage}{0.25\textwidth}
\begin{tikzpicture}[thick,scale=0.74, every node/.style={scale=0.9}]
  \begin{loglogaxis}[
    width=1.4\textwidth,
    height=1.2\textwidth,
    title style={font=\tiny},every axis title/.style={above left,at={(1,1)},draw=black,fill=white},
    title={$Q=6$},
    xlabel={DoFs},
    legend pos={south west},
    legend cell align={left},
    cycle list name=colorGPL,
    grid,
    semithick,
    ymin={5e-3},ymax={2},
    mark options={solid}
    ]

    \addplot+[solid, blue,every mark/.append style={fill=blue!80!black},mark=o] table [x=dofs,y=t_3] {data/parameters-all/parameters-q/data_2.tex};
    \addplot+[solid, blue,every mark/.append style={fill=blue!80!black},mark=x] table [x=dofs,y=t_4] {data/parameters-all/parameters-q/data_2.tex};
    \addplot+[solid, orange,every mark/.append style={fill=orange!80!black},mark=square] table [x=dofs,y=t_5] {data/parameters-all/parameters-q/data_2.tex};
    \addplot+[solid, orange,every mark/.append style={fill=orange!80!black},mark=square*] table [x=dofs,y=t_6] {data/parameters-all/parameters-q/data_2.tex};

  \end{loglogaxis}
\end{tikzpicture}
\end{minipage}
\hspace{-0.5cm}
\begin{minipage}{0.25\textwidth}
\begin{tikzpicture}[thick,scale=0.74, every node/.style={scale=0.9}]
  \begin{loglogaxis}[
    width=1.4\textwidth,
    height=1.2\textwidth,
    title style={font=\tiny},every axis title/.style={above left,at={(1,1)},draw=black,fill=white},
    title={$Q=8$},
    xlabel={DoFs},
    legend pos={south west},
    legend cell align={left},
    cycle list name=colorGPL,
    grid,
    semithick,
    ymin={5e-3},ymax={2},
    mark options={solid}
    ]

    \addplot+[solid, blue,every mark/.append style={fill=blue!80!black},mark=o] table [x=dofs,y=t_3] {data/parameters-all/parameters-q/data_3.tex};
    \addplot+[solid, blue,every mark/.append style={fill=blue!80!black},mark=x] table [x=dofs,y=t_4] {data/parameters-all/parameters-q/data_3.tex};
    \addplot+[solid, orange,every mark/.append style={fill=orange!80!black},mark=square] table [x=dofs,y=t_5] {data/parameters-all/parameters-q/data_3.tex};
    \addplot+[solid, orange,every mark/.append style={fill=orange!80!black},mark=square*] table [x=dofs,y=t_6] {data/parameters-all/parameters-q/data_3.tex};

  \end{loglogaxis}
\end{tikzpicture}
\end{minipage}

\caption{Comparison of throughputs of complex IRK and stage-parallel IRK with PRESB
or GMG as preconditioner for $k=1$ with 768
processes (16 compute nodes).}\label{fig:complex:presbvsgmg}

\end{figure}
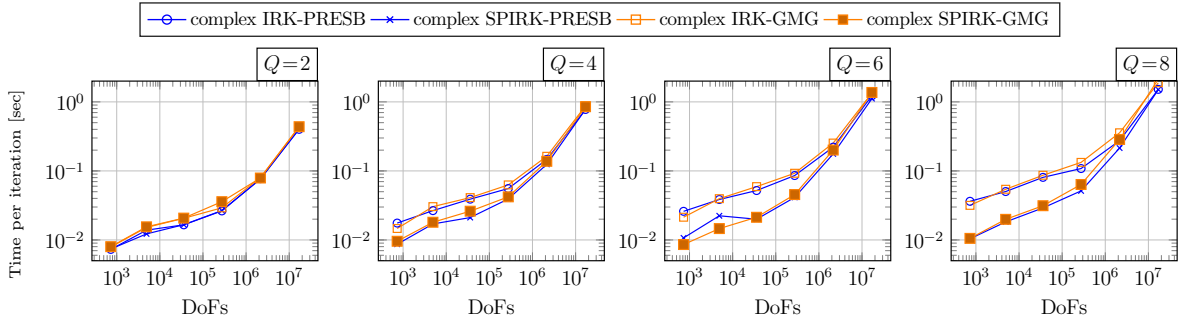

{\color{\myred}
According to our performance model, the total execution time at the scaling limit will be
dominated by the maximum number of accumulated GMG iterations each process group
has to execute. Hence, non-complex stage-parallel IRK might be advantageous for low values of $Q$,
but complex stage-parallel IRK might be competitive for high $Q$, since the number of iterations
per stage is independent of $Q$.
The results in Figure~\ref{fig:complex:noncomplexvscomplex} verify our
expectations, by comparing the times of the real-valued IRK studied before against the complex-valued IRK and the
stage-parallel IRK solvers. Surprisingly, the two stage-parallel IRK algorithms show similar times for small
problem sizes, demonstrating that the growing iteration counts of the outer solver of the
non-complex solver---indicated by the growing gap between the IRK implementations---can be
compensated by its better parallel behavior. This can also be seen in Table~\ref{tab:complex:vs}, which
shows a similar number of V-cycles
to be executed in the context of non-complex and complex stage-parallel IRK, respectively.
However, we would like to note that it is hard to
make a fair comparison, since the tolerances have different meanings
in the complex and the non-complex case. Furthermore, it was easy to integrate
efficient complex arithmetic in our small benchmarks, but this might not
be the case when the time steppers should be applied in a black-box fashion.
Against this background, a more detailed analysis for when to favor one method over the other is subject of future work.

{Figure~\ref{fig:complex:presbvsgmg} compares the solution times of PRESB and GMG
as preconditioners. Overall, PRESB seems to be superior compared to GMG 
due to fewer GMRES iterations
(see Table~\ref{tab:complex:vs}). However, the scalability is comparable. In fact,
Table~\ref{tab:complex:vs} shows that GMG needs less V-cycles in total, indicating
potential for better scalability of GMG once the
optimizations discussed in Section~\ref{sec:performance:batched} have been realized.}


Recently, Southworth et al.~\cite{southworth2022fast} presented a novel solver for IRK. They also exploit
the complex-conjugate property of the eigenvalues and solve the pairs
together. For preconditioning the blocks, they apply 2 V-cycles of (algebraic) multigrid as well.
However, it is a bit more involved to solve the blocks, since they have the form $(\Re(\lambda_i) \mathbb{I}_n - \delta t M^{-1} K)^2 + \Im(\lambda_i) \mathbb{I}_n$. The algorithm proposed in that study is,
however, inherently stage-serial as it requires the sequential solution of stage (pairs) by construction so that the algorithms presented in our publication are
not applicable there.}

%
%
%
%
%
%
%

\section{Conclusions \& outlook}\label{sec:outlook}

For distributed memory computing platforms, we have presented  {implementations of implicit Runge--Kutta algorithms, including the novel preconditioner proposed in Axelsson and Neytcheva~\cite{axelsson2020numerical}. The algorithms allow to run the
matrix-vector multiplication and preconditioner in parallel by process groups associated
with each stage. Upon a basis change---involving all stages---inner block solvers can be applied independently in a black-box fashion.
We have identified that the tensor operations
$\textbf{v}=(\mathbb{I}_Q \otimes C) \textbf{u}$ and
$\textbf{v}=(D  \otimes \mathbb{I}_n) \textbf{u}$ are the main building blocks, and have proposed efficient parallel algorithms.

We have presented a detailed performance analysis of the stage-parallel preconditioner implementation for a time-dependent
heat problem on up to 150k processes on 3k CPU compute nodes, using
state-of-the-art matrix-free geometric multigrid solvers. Furthermore, we have compared its performance to the one of
an implementation not using stage parallelism. We observed that the stage-parallel
implementation is able to significantly shift the scaling limit and reaches speedups $\le Q$. In absolute numbers, the proposed solvers and implementations make it possible to obtain high solver efficiencies down to less than 0.05 seconds per time step for four-stage IRK schemes. However, far
from the scaling limit, possible load imbalances between the solvers
of the stages and communication overhead lead to a slight drop in performance in comparison to the
non-stage-parallel IRK implementation with 13\% lower throughput on average.

{The algorithms are also applicable to the case when the system matrix
arising from implicit Runge--Kutta method is directly factorized, requiring
complex arithmetic or the solution of two-by-two blocks. This limits the
scalability in comparison to the case of the stage-parallel preconditioner. 
Results, however, show
that the lower number of iterations and the scalability balance each other,
leading to similar minimum times to solution of stage-parallel direct factorization and preconditioning for a high number of stages ($Q\le8$).

We have also discussed batching the operations of all stages instead of assigning
each stage to a process group, showing the challenge in terms of black-box interfaces.}
In future work, we plan to {study stage parallelism 
for non-linear partial differential equations,
e.g., Navier--Stokes equations. Furthermore,
we intend to improve the batched implementation as well as make further investigations
of its performance and usability within a library context, which might
require code generation~\cite{farrell2021irksome}.}


\section*{Acknowledgments}

The authors acknowledge collaboration with the \texttt{deal.II} community.
This work was supported by the Bayerisches Kompetenznetzwerk f\"ur 
Technisch-Wissenschaftliches Hoch- und H\"ochstleistungsrechnen 
(KONWIHR) through the project
``High-order matrix-free finite element implementations with
hybrid parallelization and improved data locality''. The authors
gratefully acknowledge the Gauss Centre for Supercomputing e.V.
(\url{www.gauss-centre.eu}) for funding this project by providing
computing time on the GCS Supercomputer SuperMUC-NG at
Leibniz Supercomputing Centre (LRZ, \url{www.lrz.de}) through
project id pr83te.

%
%
%

\bibliographystyle{siamplain} 
\bibliography{main}
\end{document}